\documentclass[12pt]{amsart}
\usepackage[dvips]{graphics}
\makeatletter
\theoremstyle{plain}
\newtheorem{thm}{Theorem}[section]
\numberwithin{equation}{section} 
\numberwithin{figure}{section} 
\theoremstyle{plain}
\newtheorem*{thm*}{Theorem}
\theoremstyle{plain}
\newtheorem{lem}[thm]{Lemma} 
\theoremstyle{plain}
\newtheorem{prop}[thm]{Proposition} 
\theoremstyle{plain}
\newtheorem*{prop*}{Proposition}
\theoremstyle{definition}
\newtheorem{defn}[thm]{Definition}
\theoremstyle{definition}
\newtheorem*{defn*}{Definition}
\theoremstyle{definition}
\newtheorem{example}[thm]{Example}
\theoremstyle{definition}
\newtheorem*{example*}{Example}
\theoremstyle{remark}
\newtheorem{rem}[thm]{Remark}
\theoremstyle{remark}
\newtheorem*{rem*}{Remark}

\usepackage{amsmath}
\usepackage{amsfonts}
\usepackage{amssymb}
\def\C{\mathbb{C}}
\def\Z{\mathbb{Z}}
\def\R{\mathbb{R}}
\def\P{\mathbb{P}}
\def\Q{\mathbb{Q}}
\def\I{\mathbb{I}}
\def\c{\mathcal{C}}
\def\z{\mathcal{Z}}
\def\de{\mathcal{D}}
\def\s{\square}
\def\d{\partial}
\def\D{\partial_\mathcal{B}}
\def\G{\Gamma}
\def\e{\epsilon}
\def\w{\omega}
\def\b{\bullet}
\def\X{\mathcal{X}}
\def\N{\mathcal{N}}
\def\v{\wedge}
\def\dlog{\text{dlog}}
\def\di{\text{d}}
\def\i{\iota}
\def\j{\text{\rm i}}
\def\m{\setminus}
\def\DI{{}^{\prime}\mathcal{D}}
\def\O{\Omega}
\def\r{\mathcal{R}}

\def\Cone{\text{Cone}}
\def\W{\mathcal{W}}
\def\Y{\mathcal{Y}}
\def\o{\mathcal{O}}
\def\A{\mathbb{A}}
\def\f{\textbf{f}}
\makeatother

\begin{document}

\title{The Abel-Jacobi Map for Higher Chow Groups}

\author{Matt Kerr, James D. Lewis and Stefan M\"uller-Stach}
\address{Department of Mathematics\\ UCLA, Los Angeles, CA\\ 
USA 90024-1555}
\email{matkerr@math.ucla.edu}
\address{Department of Mathematics\\ University of Alberta\\ Edmonton, 
Alberta\\ Canada  T6G 2G1}
\email{lewisjd@gpu.srv.ualberta.ca}
\address{Fachbereich Mathematik und Informatik \\ Universit\"at Mainz \\
55099 Mainz \\ Germany}
\email{mueller-stach@mathematik.uni-mainz.de}

\thanks{Second author partially supported
by a grant from the Natural Sciences and Engineering Research Council
of Canada. Third author supported
by a startup grant from McMaster University.}
\subjclass{14C25, 14C30, 14C35, 19E15}

\keywords{Abel-Jacobi map, regulator, Deligne cohomology, higher Chow group}

\begin{abstract} We construct a map between
Bloch's higher Chow groups and Deligne homology for smooth, complex 
quasiprojective varieties on the level of complexes. For complex 
projective varieties this results in a formula which generalizes
at the same time the classical Griffiths Abel-Jacobi map and the
Borel/Beilinson/Goncharov regulator type maps. \end{abstract}

\maketitle
\begin{verse}
\centerline{\it In memory of Fabio Bardelli.}
\end{verse}
\section{\textbf{Introduction}}

Let \( CH^{\bullet }(X,\bullet ) \) be the higher Chow groups as 
introduced by Bloch \cite{Bl1}, and let
\( H_{\de }^{\bullet }(X,\Z (\bullet )) \) be
Deligne cohomology. Bloch \cite{Bl2} constructed, for
\( X \) smooth, a cycle-class map
\[
c_{p,n}:\, CH^{p}(X,n)\to H_{\de }^{2p-n}(X,\Z (p)).\]

A somewhat different, but equivalent, approach using extension 
classes, is presented in \cite{D-S} and \cite{Sch}. The purpose of 
this paper is to
give an explicit description of this map in terms of currents.

More specifically, we are interested in the case where \( X \) is a 
smooth projective variety defined over \( \C \), and the higher cycle 
group in question is that of the nullhomologous cycles \( CH_{\text 
{hom}}^{p}(X,n) \). The results of this paper
pertain to this restricted setting.

As suggested by the title of this paper, we
are generalizing the classical Abel-Jacobi map
involving a membrane integral (Griffiths' prescription)
and the Borel/Beilinson/Goncharov regulator type maps
involving multiple logarithms,
to higher Chow groups. To state this 
more precisely,
we recall that there is a short exact sequence:
\[
0\to \frac{H^{2p-n-1}(X,\C )}{F^{p}H^{2p-n-1}(X,\C )+H^{2p-n-1}(X,\Z 
(p))}\to H_{\de }^{2p-n}(X,\Z (p))\]
\[\to H^{2p-n}(X,\Z (p))\bigcap F^{p}\to 0.\]
Put \[
CH_{\text {hom}}^{p}(X,n)=\ker \left\{ CH^{p}(X,n)\to H^{2p-n}_{\de 
}(X,\Z (p))\to H^{2p-n}(X,\Z (p))\right\} .\]

\begin{defn}
The induced map
\[
\Phi _{p,n}:\, CH^{p}_{\text {hom}}(X,n)\to \frac{H^{2p-n-1}(X,\C 
)}{F^{p}H^{2p-n-1}(X,\C )+H^{2p-n-1}(X,\Z (p))},\]
is called the Abel-Jacobi map.
\end{defn}
Let \( m=\dim X \). By Poincar\'e and Serre duality, we will think of 
this map in the form:
\[
\Phi _{p,n}:\, CH^{p}_{\text {hom}}(X,n)\to \frac{\left\{ 
F^{m-p+1}H^{2m-2p+n+1}(X,\C )\right\} ^{\vee }}{H_{2m-2p+n+1}(X,\Z 
(m-p))}.\]
We will use the cubical complex description of \( CH^{p}(X,n) \) 
throughout this paper (see \( \S 2). \) Cycles live in \( X\times \s 
^{n} \), where \( \s ^{n}:=(\P _{\C }^{1}\m \{1\})^{n} \)
has coordinates \( (z_{1},\ldots ,z_{n}) \), and there are projection 
maps
\( \pi _{X}:\, X\times \s ^{n}\to X \) and \( \pi _{\s }:\, X\times 
\s ^{n}\to \s ^{n} \).
Consider an irreducible subvariety \( \z \subset X\times \s ^{n} \), 
of codimension
\( p \), and a form \( \omega \in F^{m-p+1}\Omega _{X^{\infty 
}}^{2m-2p+n+1}(X) \).

One considers the current associated to \( \z \), defined by its 
action on
\( \omega \):\footnote{Strictly speaking, in terms of homology
the Tate twist should be \(\frac{1}{(2\pi \j)^{m-p+n}}\).
However we felt it was more natural to give a cohomological formulation of the
map, where the extra \((2\pi \j)^{m}\) is lost in passage 
from Deligne homology to cohomology, as indicated in 5.6.}
\[
\frac{1}{(2\pi \j)^{n-p}}\left[ \int_{{}_{\z \m \{\z \cap \pi 
^{-1}_{\s }([-\infty ,0]\times \s ^{n-1})\}}}\pi ^{*}_{\s }\left( 
(\log z_{1})\dlog z_{2}\v \cdots \v \dlog z_{n}\right) \v \pi 
_{X}^{*}\w \right. \]
\[
+\, (-2\pi {\j})\int_{\begin{array}{clcr}
{}_{\{\z \cap \pi _{\s }^{-1}[-\infty ,0]\times \s ^{n-1}\}}\\ {}^{\m 
\{\z \cap \pi _{\s }^{-1}([-\infty ,0]^{2}\times \s 
^{n-2})\}}\end{array}}\pi _{\s }^{*}\left( (\log z_{2})\dlog z_{3}\v 
\cdots \v \dlog z_{n}\right) \v \pi _{X}^{*}\w \]
\[
+\cdots +\, (-2\pi {\j})^{n-1}
\int_{\begin{array}{clcr}
{}_{\{\z \cap \pi ^{-1}_{\s }([-\infty ,0]^{n-1}\times \s )\}}\\ 
{}^{\m \{\z \cap \pi_{\s }^{-1}([-\infty ,0]^{n})\}}\end{array}}\pi 
_{\s }^{*}(\log z_{n})\v \pi _{X}^{*}\w \]
\[
\left. +\, (-2\pi {\j})^{n}\int _{\zeta }\pi _{X}^{*}\w \right] ,\]
where the latter term is a membrane integral, and \( \log z_{i} \) 
represents the branch of the logarithm with argument in 
\((-\pi,\pi)\) (same for every \( z_{i} \)).

\begin{thm} \( \Phi _{p,n} \) is induced by the above current.
\end{thm}

The plan of the paper is as follows. After reviewing the classical 
situation (\( n=0 \)), we arrive at the formula for the \( AJ \) map, 
based on a cup-product calculation at the generic point. This was the 
point of view adopted by the second author. Around the same time the 
first author arrived at the formula based on a morphism of complexes. 
This is fully explained in \cite{Ke}, and the relevant points are 
explained here. By comparing the extension class \cite{Sch}
construction with the above formula, we arrive at
the above theorem. We are grateful to S. Bloch for remarking that the 
\( AJ \) map can be described in terms of the dilogarithm associated 
to the ``Totaro'' cycles that are discussed in \cite{Bl3}. This led 
to the inclusion of an example in \( \S 5.7 \) below.

\section{\textbf{Some definitions}}

\subsection{Higher Chow groups} (\cite{Bl1})

Let \( W/\C \) be a quasiprojective variety. Put \( Z^{k}(W)= \)free 
abelian
group generated by subvarieties of codimension \( k \) in \( W \). 
Consider
the \( n \)-simplex:
\[
\Delta ^{n}=\text {Spec}\left\{ \frac{\C [t_{0},\ldots 
,t_{n}]}{(1-\sum _{j=0}^{n}t_{j})}\right\} \simeq \C ^{n}.\]
We set \[
Z^{k}(W,n):=\left\{ \xi \in Z^{k}(W\times \Delta ^{n})\, \, \right| \]
\[
\left. \text {every\, component\, of\, }\xi \text {\, meets\, all\, 
faces\, }\{t_{i_{1}}=\ldots =t_{i_{\ell }}=0,\, \ell \geq 1\}\text 
{\, properly}\right\} .\]
Note that \( Z^{k}(X,0)=Z^{k}(W) \). Now set \( \partial _{j}:\, 
Z^{k}(W,n)\to Z^{k}(W,n-1) \)
to be the restriction to the \( j^{\text {th}} \) face (given by \( 
t_{j}=0 \)).
The boundary map \[
\delta =\sum _{j=0}^{n}(-1)^{j}\partial _{j}\, :\, Z^{k}(W,n)\to 
Z^{k}(W,n-1),\]
satisfies \( \delta ^{2}=0 \).

\begin{defn}
\cite{Bl1} \( CH^{k}(W,\bullet )= \) homology of \( \left\{ 
Z^{k}(W,\bullet ),\delta \right\} \).
\end{defn}

We put \( CH^{k}(W):=CH^{k}(W,0). \)

\subsubsection*{Cubical Version}

Let \( \s ^{n}:=(\P ^{1}\m \{1\})^{n} \) with coordinates \( z_{i} \) 
and
\( 2^{n} \) codimension-one faces obtained by setting \( 
z_{i}=0,\infty \),
and boundary maps \( \partial =\sum (-1)^{i-1}(\partial 
^{0}_{i}-\partial ^{\infty }_{i}) \),
where \( \partial _{i}^{o} \), \( \partial _{i}^{\infty } \) denote 
the restriction
maps to the faces \( z_{i}=0 \), \( z_{i}=\infty \) respectively. The 
rest of the definition is completely analogous except that one has to 
divide out degenerate cycles. The precise description is given
in 5.2 in \S 5. It is known that both complexes are quasiisomorphic.

\subsection{Deligne cohomology}

Working in the analytic topology, we introduce the Deligne complex 
(for any
subring \( \A \subseteq \R \))
\[
\A _{\de }(k):\, \, \A (k)\to \begin{array}[t]{c}
\underbrace{\o _{X}\to \O ^{1}_{X}\to \cdots \to \O ^{k-1}_{X}}\\
^{\text {call\, this\, }\O _{X}^{\bullet <k}}
\end{array}.\]

\begin{defn}
Deligne cohomology is given by the hypercohomology:
\[
H^{i}_{\de }(X,\A (k))=\mathbb {H}^{i}(\A _{\de }(k)).\]
\end{defn}
From the short exact sequence
\[
0\to \O _{X}^{\bullet <k}[-1]\to \A _{\de }(k)\to \A (k)\to 0,\]
one has the short exact sequence
\[
0\to \frac{H^{i-1}(X,\C )}{H^{i-1}(X,\A (k))+F^{k}H^{i-1}(X,\C )}\to 
H_{\de }^{i}(X,\A (k))
\]
\[
\to H^{i}(X,\A (k))\cap F^{k}H^{i}(X,\C )\to 0.\]
In particular, the case \( (\A ,i,k)=(\Z ,2p-n,p) \) gives the exact 
sequence
in the introduction.

\subsection{Deligne homology}

In this part, we follow \cite{Ja} rather closely. Let \( \mu :\, 
A^{\bullet }\to B^{\bullet } \)
be a morphism of complexes. Then the cone complex is given by
\[
\text {Cone}(A^{\bullet }\buildrel \mu \over \to B^{\bullet 
}):=A^{\bullet }[1]\oplus B^{\bullet },\]
where the differential \( \delta :\, A^{\bullet +1}\oplus B^{\bullet 
}\to A^{\bullet +2}\oplus B^{\bullet } \)
is given by \( \delta (a,b)=(-da,\mu (a)+db). \)\\

We introduce some notation: \( \O ^{p,q}_{X^{\infty }}= \) sheaf of 
\( C^{\infty } \)
\( (p,q) \)-forms on \( X \); \( '\O _{X^{\infty }}^{p,q}= \) sheaf 
of distributions
over \( \O _{X^{\infty }}^{-p,-q} \). Thus for an open set \( 
U\subset X \),
an element of \( '\O _{X^{\infty }}^{p,q}(U) \) is a continuous 
linear functional
on the compactly supported forms \( \Gamma _{c}(U,\O _{X^{\infty 
}}^{-p,-q}) \).

\subsubsection*{Key Example I:}

Any \( C^{\infty } \) \( (p,q) \)-form \( \eta \) gives a section of 
\( '\O _{X^{\infty }}^{p-m,q-m} \)
by the formula:
\[
\omega \mapsto l(\eta )(\omega )=\frac{1}{(2\pi {\j})^{m}}\int 
_{X}\eta \v \omega .\]

\subsubsection*{Key Example II:}

Any piecewise smooth oriented topological \( r \)-chain \( \xi \) on 
\( X \)
gives a current \( \delta _{\xi } \) in \( \oplus_{p+q=-r}{}^{\prime}
\Omega_{X^{\infty }}^{p,q}(X) \)
by the formula
\[
\omega \mapsto \epsilon (\xi )(\omega )=\int _{\xi }\omega .\]

\( '\Omega _{X^{\infty }}^{\bullet ,\bullet } \) and \( \Omega 
_{X^{\infty }}^{\bullet ,\bullet } \)
naturally form double complexes, where if \( D \) is a current acting 
on an
\( r \)-form, and given an \( (r-1) \)-form \( \omega \), then \( 
dD(\omega )=(-1)^{r-1}D(d\omega ) \),
\( d=\partial +\bar{\partial } \). [Warning: The definition
of \(dD\) differs from \cite{Ja} by a minus sign.] 
One has Hodge filtrations:
\[
F^{i}\Omega _{X^{\infty }}^{\bullet }=\bigoplus _{p+q=\bullet ,\, 
p\geq i}\Omega ^{p,q}_{X^{\infty }}\, ,\, \, \, \, \, \, \, \, \, \, 
\, \, \, F^{i}{}'\O _{X^{\infty }}^{\bullet }=\bigoplus _{p+q=\bullet 
,\, p\geq i}{}'\O ^{p,q}_{X^{\infty }}\, .\]

Let \( (C_{\bullet }(X,\Z (k)),d) \) be the complex of singular \( 
C^{\infty } \)-chains
with coefficients in \( \Z (k) \), and put \( 'C^{i}=C_{-i} \) 
{[}with differential
\( (-1)^{i+1}d:\, C_{-i}\to C_{-i-1} \){]}. One has a morphism of 
complexes
\[
\epsilon :\, {}'C^{\bullet }(X,\Z (k))\to {}'\Omega ^{\bullet 
}_{X^{\infty }}(X).\]
Put \[
M_{\de }^{\bullet }=\text {Cone}\left\{ {}'C^{\bullet }(X,\Z 
(p-m))\oplus F^{p-m}{}'\O ^{\bullet }_{X^{\infty }}(X)\buildrel 
{\epsilon -l}\over \to {}'\O _{X^{\infty }}^{\bullet}(X)\right\} 
[-1].\]
The homology of this complex, at \( \bullet =2p-n-2m \), viz., \( 
'H_{\de }^{2p-n-2m}(X,\Z (p-m)) \),
is precisely the Deligne homology:
\[
H^{\de }_{2m-2p+n}(X,\Z (m-p)):={}'H_{\de }^{2p-n-2m}(X,\Z 
(p-m))\simeq H_{\de }^{2p-n}(X,\Z (p)).\]
(Poincar\'e duality.)

\begin{rem}
A class in \( H^{\de }_{2m-2p+n}(X,\Z (m-p)) \) is represented by a 
triple
\[
(a,b,c)\, \in \, {}'C^{2p-n-2m}(X,\Z (p-m))\bigoplus F^{p-m}{}'\O 
^{2p-n-2m}_{X^{\infty }}(X)\bigoplus {}'\O _{X^{\infty 
}}^{2p-n-2m-1}(X),\]
where \( da=0 \), \( db=0 \), and \( a-b+dc=0 \). Via Poincar\'e 
duality,
this corresponds to \( [a]\in H^{2p-n}(X,\Z (p)) \), \( [b]\in 
F^{p}H^{2p-n}_{\text {DR}}(X,\C ) \),
with \( [a]=[b] \) in \( H^{2p-n}(X,\C ) \). Now suppose that \( 
[a]=[b]=0 \).
Then \( a=da_{0} \), and from Hodge theory, \( b=db_{0} \), where \( 
b_{0}\in F^{p-m}{}'\O _{X^{\infty }}^{2p-n-2m-1}(X) \).
Thus \( d(a_{0}-b_{0}+c)=0 \), and \( [a_{0}-b_{0}+c] \) represents 
the corresponding
class in
\[
\frac{H^{2p-n-1}(X,\C )}{F^{p}H^{2p-n-1}(X,\C )+H^{2p-n-1}(X,\Z 
(p))}\simeq \frac{F^{m-p+1}H^{2m+n-2p+1}(X,\C )^{\vee 
}}{H_{2m+n-2p+1}(X,\Z (m-p))}.\]
By Hodge-type considerations, the action of the current \( b_{0} \) on
\[
F^{m-p+1}H^{2m+n-2p+1}(X,\C )\]
is zero. Thus the action of the closed current \[
a_{0}-b_{0}+c\text {\, \, on\, \, }F^{m-p+1}H^{2m+n-2p+1}(X,\C ),\]
is the same as the action of \( a_{0}+c \).
\end{rem}

We also need a slightly expanded version of Deligne homology for the 
smooth
quasiprojective case. Let \( Z \) be a smooth quasiprojective with 
good compactification
\( \bar{Z} \) (with NCD \( E \)). \( C_{\bullet }(\bar{Z},\A (k)) \) 
is the
complex of singular \( C^{\infty } \)-chains in \( \bar{Z} \) with 
coefficients
in \( \A (k) \), and \( 'C^{i}=C_{-i} \). Let \[
'C^{\bullet }(\bar{Z},E,\A (k))='C^{\bullet }(\bar{Z},\A (k))\left/ 
{}'C_{E}^{\bullet }(\bar{Z},\A (k))\right. \]
where \( 'C^{\bullet }_{E}(\bar{Z},\A (k))\subset {}'C^{\bullet 
}(\bar{Z},\A (k)) \)
is the subcomplex of chains supported on \( E \).

Deligne homology \( 'H_{\de }^{\bullet }(Z,\A (k)) \), as defined
in \cite{Ja}, is given by the cohomology
of the complex
\[
\text {Cone}\left( {}'C^{\bullet }(\bar{Z},E,\A (k))\oplus F^{k}{}'\O 
^{\bullet }_{\bar{Z}^{\infty }}\left\langle E\right\rangle 
(\bar{Z})\buildrel {\epsilon -l}\over \to {}'\O _{\bar{Z}^{\infty 
}}^{\bullet }\left\langle E\right\rangle (\bar{Z})\right) [-1],\]
where \( \epsilon \) and \( l \) are the natural maps of complexes.
(The precise description of \( \epsilon \) is given
in \cite{Ja}, and the required foundational material
can be found in \cite{Ki}.) Here
we define \( \O ^{\bullet }_{\bar{Z}}\left\langle 
E\right\rangle =\O^{\bullet }_{\bar{Z}}(\log E) \)
to be the de Rham complex of meromorphic forms on \( \bar{Z} \), 
holomorphic
on \( U=\bar{Z}-E \), with at most logarithmic poles along \( E \).
Also, \( \O _{\bar{Z}^{\infty }}^{\bullet }\left\langle 
E\right\rangle =\O _{\bar{Z}}^{\bullet }\left\langle E\right\rangle 
\otimes _{\O _{\bar{Z}}^{\bullet }}\O _{\bar{Z}^{\infty }}^{\bullet } 
\),
and \( '\O^n_{\bar{Z}^{\infty}}\left\langle E\right\rangle \) is defined by the 
equivalent sheaves
\[
\left. '\O^n_{\bar{Z}^{\infty}} \right/ '\O^n_{\bar{Z}^{\infty}}(on E)
\cong \mathcal{D} \left( \O^{-n}_{\bar{Z}^{\infty}}(null E) \right)
\cong \bigoplus_{p+q=n} \O^{p+m}_{\bar{Z}}\left\langle E\right\rangle
\otimes_{\mathcal{O}_Z} '\O^{0,q}_{\bar{Z}^{\infty}}. \]
There is thus a map of complexes \( '\O^{\b}_{\bar{Z}^{\infty}} \to
'\O^{\b}_{\bar{Z}^{\infty}}\left\langle E\right\rangle \) which is 
surjective at each term.
\( '\O _{\bar{Z}^{\infty }}^{\bullet }\left\langle E\right\rangle 
=\O _{\bar{Z}}^{\bullet }\left\langle E\right\rangle \otimes _{\O 
_{\bar{Z}}^{\bullet }}{}'\O _{\bar{Z}^{\infty }}^{\bullet } \).
The corresponding Hodge filtrations are 
\( F^{i}\O _{\bar{Z}^{\infty }}^{\bullet 
}\left\langle E\right\rangle =\left\{ F^{i}\O _{\bar{Z}}^{\bullet 
}\left\langle E\right\rangle \right\} \otimes _{\O 
_{\bar{Z}}^{\bullet }}\O _{\bar{Z}^{\infty }}^{\bullet } \),
and \( F^{i}{}'\O _{\bar{Z}^{\infty }}^{\bullet }\left\langle 
E\right\rangle =\left\{ F^{i+m}\O _{\bar{Z}}^{\bullet }\left\langle 
E\right\rangle \right\} \otimes _{\mathcal{O}_{\bar{Z}}}{}'\O 
_{\bar{Z}^{\infty }}^{0, \bullet } \).
As is well-known (see \cite{Ja}), there are filtered quasiisomorphisms

\[
\left( \O _{\bar{Z}}^{\bullet }\left\langle E\right\rangle 
,F^{i}\right) \hookrightarrow \left( \O ^{\bullet }_{\bar{Z}^{\infty 
}}\left\langle E\right\rangle ,F^{i}\right) \hookrightarrow \left( 
{}'\O _{\bar{Z}^{\infty }}^{\bullet }\left\langle E\right\rangle 
[-2m],F^{i-m}\right) .\]

\section{\textbf{Review of the classical situation (\protect\(
n=0\protect \))}}

General references for this section are \cite{E-V} and \cite{Ja}.
For a codimension \( p \) cycle \( \z \) on \( X \), there is the 
localization sequence of mixed Hodge structures

\begin{equation}
0\to H^{2p-1}(X,\Z(p)){\buildrel\beta\over\to} H^{2p-1}(X\m 
|\z|,\Z(p)) \to H^{2p}_{|\z|}(X,\Z(p))\to H^{2p}(X,\Z(p))\\
\end{equation}
The map \( \beta \) induces the isomorphism:

\begin{equation}
\frac{H^{2p-1}(X,\C)}{F^pH^{2p-1}(X,\C)}\simeq\frac{H^{2p-1}(X\m 
|\z|,\C)}{F^pH^{2p-1}(X\m |\z|,\C)} ,\\
\end{equation}\\
and hence the isomorphism

\begin{equation}
J^p(X)\simeq 
\frac{H^{2p-1}(X\m|\z|,\C)}{H^{2p-1}(X,\Z(p))+F^pH^{2p-1}(X\m|\z|,\C)} 
.\\
\end{equation}
Next, for \( \z \in Z^{p}_{\text {hom}}(X), \) the fundamental class 
\( c_{\Z }(\z ) \)
is the image of a class \( \widetilde{c_{\Z }(\z )}\in H^{2p-1}(X\m 
|\z |,\Z (p)) \),
uniquely determined up to \( \text {Im}(\beta ) \). Since \( 
\widetilde{c_{\Z }(\z )} \)
defines a class in \( H^{2p-1}(X\m |\z |,\C ) \) (still denoted by \( 
\widetilde{c_{\Z }(\z )} \)),
we end up with a corresponding class \( \Psi _{p}(\z )\in J^{p}(X) \) 
via the
isomorphisms above. We use this as our initial

\begin{defn}
\( \Psi _{p}:\, Z_{\text {hom}}^{p}(X)\to J^{p}(X) \) is called the 
Abel-Jacobi
map.
\end{defn}

\subsection{Comparison to Carlson's Abel-Jacobi map}

The exact sequence \( (3.1) \) yields an extension
\[
0\to H^{2p-1}(X,\Z (p))\to \textbf {E}\to \Z (0)\to 0\]
via pullback, where \( \mathbf{E} \) is abstractly identified with \( 
H^{2p-1}(X,\Z (p))\oplus \Z \widetilde{c_{\Z }(\z )} \),
equipped with an (integral) retraction \( r_{\Z }:\, \mathbf{E}\to 
H^{2p-1}(X,\Z (p)) \)
killing \( \widetilde{c_{\Z }(\z )} \). This \( r_{\Z } \) extends to 
a map \( \mathbf{E}_{\C }\to H^{2p-1}(X,\C ) \).

If \( \widetilde{c_{F}(\z )}\in E_{\C }\subseteq F^{p}H^{2p-1}(X\m 
|\z |,\C ) \)
is another lift of the fundamental class of \( \z \) respecting the
Hodge filtration, then Carlson's prescription \cite{C} is \( r_{\Z 
}(\widetilde{c_{F}(\z )})\in J^{p}(X) \),
its image under the retraction. Since \( \widetilde{c_{\Z }(\z )} \) 
and \( \widetilde{c_{F}(\z )} \)
both lift the fundamental class, their difference lifts to an element 
\( \xi \in H^{2p-1}(X,\C ) \). We write this \( \widetilde{c_{F}(\z 
)}=\widetilde{c_{\Z }(\z )}+\xi \);
applying \( r_{\Z } \) shows that \( \xi =r_{\Z }(\widetilde{c_{F}(\z 
)}) \),
and hence that \( \widetilde{c_{F}(\z )}=\widetilde{c_{\Z }(\z 
)}+r_{\Z }(\widetilde{c_{F}(\z )}) \).
Thus
\[
\widetilde{c_{\Z }(\z )}+r_{\Z }(\widetilde{c_{F}(\z )})\equiv 0\text 
{\, \, \, modulo\, \, }F^{p}H^{2p-1}(X\m |\z |,\C ).\]
So we have the

\begin{prop}
Carlson's Abel-Jacobi map is the same as \( \Psi _{p} \), up to sign.
\end{prop}

\subsection{Comparison to classical AJ map}

We proceed by comparing \( \Psi _{p} \) with the Deligne cycle-class 
map in
Prop. 3.3, and then identifying the latter with the classical AJ in 
Prop. 3.4.

Recall the diagram of exact sequences
\[
\begin{array}{ccccccccc}
0&\to&Z^{p}_{\text{hom}}(X)&\to&Z^{p}(X)&\to&
Z^{p}(X)/Z^{p}_{\text{hom}}(X)&\to&0 \\
&\\
&&\quad \downarrow {\Phi_{p,0}}&&
\quad \downarrow {cl_{p,0}}&&\downarrow\\
&\\
0 &\to & J^p(X) &\to & H^{2p}_{\de}(X,\Z(p))& \to& H^{2p}(X,\Z(p)) 
\bigcap F^p & \to& 0 \end{array}
\]

\begin{prop}
\( \Psi _{p}=\Phi _{p,0}. \)
\end{prop}
\begin{proof}
See \cite{E-V}.
\end{proof}
Finally, working with Deligne homology:
\begin{prop}
\( \Phi _{p,0} \) coincides with the classical Abel-Jacobi map.
\end{prop}
\begin{proof} See \cite{Ja}.
\end{proof}

\section{\textbf{A localization argument (first construction of 
\protect\( AJ\protect \))}}

If \( \z \in X\times \s ^{n} \) is irreducible and of codimension \( 
p \),
consider \( V:=\overline{\pi _{*}(\z )}\subset X \), which we assume 
has dimension
\( m+n-p \). The \( i \)th coordinate projections \( \z \to \s \), 
determine
rational functions on \( \z \). Taking the norm of a symbol in Milnor 
\( K \)-theory,
after passing to the relevant functions fields, reduces to the 
situation of
rational functions \( \{f_{1},\ldots ,f_{n}\} \) on \( V \). Of 
course when
\( n=1 \), we are dealing with the usual norm \( N:\, \C (\z 
)^{\times }\to \C (V)^{\times } \).
On an open set \( U_{V}\subset V \), we have elements \( f_{j}\in 
H^{0}(U_{V},\mathcal{O}_{U_{V}}^{\times })=CH^{1}(U_{V},1). \)
Let \( U=X\m (V\m U_{V}) \). One has this commutative diagram:

\[
\begin{array}{ccc}
CH^1(U_V,1)^{\otimes n} & \to & H^1_{\de}(U_V ,\Z(1))^{\otimes n} \\
&\\
\quad \downarrow {\cup} & & \quad\downarrow {\cup} \\
&\\
CH^n(U_V, n) & \to & H^n_{\de} (U_V, \Z(n)) \\
&\\
\downarrow & & \downarrow \\
&\\
CH^p(U,n) & \to & H^{2p-n}_{\de}(U,\Z(n)) \\
&\\
\uparrow & & \uparrow \\
&\\
CH^p(X,n) & \to & H^{2p-n}_{\de} (X,\Z(n)) \\
\end{array} \]
We come up with a formula for the regulator at the generic point, 
based on the cup product formula in Deligne cohomology.
We need the following basic

\begin{lem}
Let \( Y \) be a smooth quasiprojective variety, \( f:\, Y\to \P ^{1} 
\) a
dominant morphism, and let \( \eta \in E_{\bar{Y}}^{2m-1}. \) Then
\[
\int _{\bar{Y}}\frac{df}{f}\v \eta \, =\, (2\pi {\j})\int 
_{f^{-1}[-\infty,0 ]}\eta \, +\, d[T_{\log f}](\eta ),\]
where \( T_{\log f}(\mu ):=\int _{\bar{Y}\m f^{-1}[-\infty,0 ]}\mu 
\log f \)
and \( [-\infty ,0]=\R ^{-} \) is oriented along \( \R ^{-} \) so 
that \( \partial [-\infty ,0]=\{0\}-\{\infty \} \).
\end{lem}
\begin{proof}
One observes that \( \log f \) is single-valued in \( \bar{Y}\m 
f^{-1}[-\infty,0 ] \),
where \( \log \) is the branch with \( \arg \in (-\pi ,\pi ) \); so 
we have
\[
d((\log f)\eta )\, =\, \frac{df}{f}\v \eta +(\log f)d\eta \]
there. Let \( B_{\epsilon } \) be an \( \epsilon \)-band angular 
sector neighborhood
of \( [-\infty ,0]=\R ^{-} \) in \( \P ^{1}\), with boundary \( 
C_{\pm }(\epsilon)\).

\includegraphics{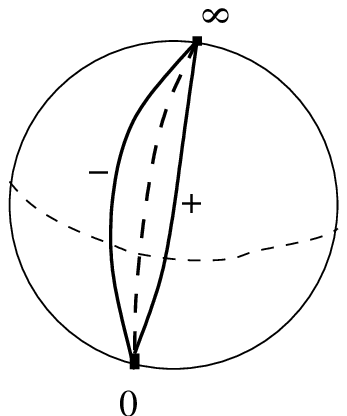}

Put \( D_{\epsilon }:=f^{-1}(B_{\epsilon }) \), and \( L_{\pm 
}(\epsilon )=\partial D_{\epsilon }=f^{-1}(C_{\pm }(\epsilon )). \)
Then
\[
\int _{\bar{Y}}\frac{df}{f}\v \eta \, =\, \int _{\bar{Y}\m 
f^{-1}[-\infty ,0]}d\left( (\log f)\eta \right) \, -\, \int 
_{\bar{Y}}(\log f)d\eta \]
\[
= \lim _{\epsilon \to 0^{+}}\int _{\bar{Y}\m D_{\epsilon }}d\left( 
(\log f)\eta \right) \, \, +\, d[T_{\log f}](\eta )\]
\[
=\, -\lim _{\epsilon \to 0^{+}}\int _{L_{\pm }(\epsilon )}(\log 
f)\eta \, \, +\, d[T_{\log f}](\eta ).\]
But \[
\lim _{\epsilon \to 0^{+}}L_{\pm }(\epsilon )\, =\, f^{-1}[-\infty 
,0]-f^{-1}[-\infty ,0]\, =0.\]
Since we pick up a period on \( \log \), we arrive at:
\[
\int _{\bar{Y}}\frac{df}{f}\v \eta \, =\, 2\pi {\j}\int 
_{f^{-1}[-\infty,0 ]}\eta \, \, +\, d[T_{\log f}](\eta ),\]
as was to be shown.
\end{proof}
We now want to consider the following setting. Let \( V \) be an 
irreducible
complex projective variety (in particular a component of the \( V \) 
from before
the lemma). Let \( f_{1},\ldots ,f_{n}\in \C (V)^{\times } \), and 
put \( D:=\bigcup ^{n}_{j=1}|(f_{j})|\bigcup V_{\text {sing}} \),
\( U_{V}:=V\m D \). Consider the pair \( (\tilde{V},\tilde{D}) \), 
where \( \tilde{V} \)
is a smooth projective variety, \( \tilde{D} \) a normal crossing 
divisor,
and \( \tilde{V}\m \tilde{D}=V\m D \). The Deligne (co)homology of \( 
U_{V} \)
can be computed in terms of the Deligne complex of the pair \( 
(\tilde{V},\tilde{D}) \).
We may assume that \( f_{1},\ldots ,f_{n}:\, \tilde{V}\to \P ^{1} \) 
are dominant
morphisms. Notice that \( \gamma _{j}:=f_{j}^{-1}[-\infty ,0] \) are 
(Borel-Moore)
cycles on \( U_{V}. \) Let \( v=\dim V\, (=m+n-p) \) and define \( 
T_{f_{j}}(\mu )=2\pi {\j}\int _{\gamma _{j}}\mu \),
\( \O _{f_{j}}(\mu )=\int _{\tilde{V}}\frac{df_{j}}{f_{j}}\v \mu \), 
\( R_{f_{j}}(\mu )=\int _{\tilde{V}}(\log f_{j})\mu . \)

Before stating our next result, we recall the multiplication table 
(\cite{E-V}), pertaining to
\[
\bigcup _{\alpha }:\, {\A}_{\de }(p)\otimes {\A}_{\de }(q)\to 
{\A}_{\de }(p+q),\]
for any given \( \alpha \in \R \). Here \( {\A}_{\de }(p) \) is 
defined in \cite{E-V} in
terms of a cone complex, which is quasiisomorphic to
the same labelled complex \( {\A}_{\de }(p) \) that we introduced in 
\S 2.

\begin{equation}
\end{equation}
\[
\begin{array}[b]{ccccccc}
& | & a_{q} & | & f_{q} & | & \w _{q}\\
-- & & ----- & & ----- & & -----\\
a_{p} & | & a_{p}\cdot a_{q} & | & 0 & | & (1-\alpha )\cdot 
a_{p}\cdot \w _{q}\\
-- & & ----- & & ----- & & -----\\
f_{p} & | & 0 & | & f_{p}\v f_{q} & | & (-1)^{\deg f_{p}}\cdot \alpha 
\cdot f_{p}\cdot \w _{q}\\
-- & & ----- & & ----- & & -----\\
\w _{p} & | & \alpha \cdot \w _{p}\cdot a_{q} & | & (1-\alpha )\cdot 
\w _{p}\v f_{q} & | & 0
\end{array}\]
Evidently, up to homotopy, \( \bigcup _{\alpha } \) is independent of 
\( \alpha \in \R \)
\cite{E-V}. Now put \( \bigcup =\bigcup _{\alpha =0} \). We have the 
following

\begin{prop}
(i) For each \( j \), the triple \( (T_{f_{j}},\O _{f_{j}},R_{f_{j}}) 
\) defines
a class
\[
\left\{ (T_{f_{j}},\O _{f_{j}},R_{f_{j}})\right\} \in H^{\de 
}_{2v-1}(U_{V},\Z (v-1))\begin{array}[b]{c}
_{\text {PD}}\\
\simeq \end{array}H_{\de }^{1}(U_{V},\Z (1)).\]
(ii) Via the cup product:
\[
\bigcup _{j=1}^{n}\left\{ (T_{f_{j}},\, \O _{f_{j}},\, 
R_{f_{j}})\right\} =\left\{ (T_{\f },\, \Omega _{\f },\, R_{\f 
})\right\} ,\]
where \( \f =(f_{1},\ldots ,f_{n}) \),
\[
T_{\f }(\mu )=(2\pi {\j})^{n}\int _{(f_{1}\times \cdots \times 
f_{n})^{-1}[-\infty,0]^{n}}\mu \, \, ,\, \, \, \, \Omega _{\f }(\mu 
)=\int _{V}\frac{df_{1}}{f_{1}}\v \cdots \v \frac{df_{n}}{f_{n}}\v 
\mu ,\]
and where
\[
R_{\f }(\omega )\, =\, \left[ \int _{V\m f_{1}^{-1}[-\infty,0]}(\log 
f_{1})\frac{df_{2}}{f_{2}}\v \cdots \v \frac{df_{n}}{f_{n}}\v \w 
\right. \]
\[
+\, (-2\pi {\j})\int _{f^{-1}_{1}[-\infty,0]\m (f_{1}\times 
f_{2})^{-1}[-\infty,0]^{2}}(\log f_{2})\frac{df_{3}}{f_{3}}\v \cdots 
\v \frac{df_{n}}{f_{n}}\v \w \]
\[
\left. +\cdots +\, (-2\pi {\j})^{n-1}\int _{(f_{1}\times \cdots 
\times f_{n-1})^{-1}[-\infty,0]^{n-1}\m (f_{1}\times \cdots \times 
f_{n})^{-1}[-\infty,0]^{n}}(\log f_{n})\w \right] .\]
\end{prop}
\begin{proof}
Part (i) is immediate from Lemma 4.1 and part (ii) uses the 
multiplication table above, the cone complex
description of Deligne homology together with Poincar\'e
duality, and induction on \( n \).
\end{proof}
\begin{rem}
As a consequence of part (ii) above, we have the Deligne homology 
relation

\begin{equation}
\Omega_{\f} = T_{\f}+ d[R_{\f}] ,\\
\end{equation}
as currents acting on forms that are compactly supported on \( U_{V} 
\). Using
induction, the proof of Lemma 4.1 can be generalized, which extends 
the above
formula to act on forms on \( V \):
\end{rem}
\begin{prop}
Consider (dominant) morphisms \( f_{1},\ldots ,f_{n} \) from \( V \) 
to \( \P ^{1} \),
in general position and put \[
R_{\partial \f }=\sum _{j=1}^{n}(-1)^{j-1}R_{\{f_{1},\ldots 
,\hat{f}_{j},\ldots ,f_{n}\}}\left|_{_{(f_{j})}}\right. .\]
Then
\[
\Omega _{\f }=T_{\f }+d[R_{\f }]\pm (2\pi {\j})R_{\partial \f }.\]
\end{prop}

\begin{example}
Suppose we are given a higher Chow cycle \( \z =\sum _{\alpha 
}(f_{\alpha },V_{\alpha })\in CH_{\text {hom}}^{p}(X,1) \).
Then by Hodge theory, we have \( \sum _{\alpha }\O _{f_{\alpha }}=dS 
\), where by Poincar\'e duality \( S\in F^{p} \), (which plays the
role of \( b_{0} \) in remark 2.3) acts as the zero current on \( 
F^{m-p+1}H^{2m-2p+2}(X) \).
Note that \( \gamma :=\sum _{\alpha }\gamma _{\alpha } \) bounds a 
chain \( \zeta , \)
thus \( T_{\z }:=\sum _{\alpha }T_{f_{\alpha }}= -2\pi\j 
d\delta_{\zeta } \). Taking the coboundary (see \S 2), viz.,
\[ \delta (-\delta_{\zeta },S,0) = (-T_{\z},-\sum_{\alpha}\O_{\alpha},
-2\pi\j \delta_{\zeta}- S),
\]
this leads us to \[(T_{\z },\, \sum _{\alpha }\O _{f_{\alpha }},\, 
\sum _{\alpha }R_{f_{\alpha }})\sim (0,\, 0,\, \sum _{\alpha 
}R_{f_{\alpha }}- 2\pi\j\delta_{\zeta }-S) \]
in Deligne homology. By applying Poincar\'e duality, this leads us to 
Levine's formula \cite{Lev1} for the regulator on \( K_{1} \), induced 
by
\[
\w \in F^{m-p+1}\Omega _{X^{\infty }}^{2m-2p+2}(X)
\]
\[\mapsto \frac{1}{(2\pi {\j})^{m-p+1}}\left( \sum _{\alpha }\int 
_{V_{\alpha }\m f_{\alpha }^{-1}[-\infty,0]}(\log f_{\alpha })\w \, 
\, \, -\, 2\pi {\j}\int _{\zeta }\w \right) .\]
[Note: In Levine's formula (op. cit.), the \( -2\pi {\j}\int _{\zeta 
}\w\) is replaced by \(+ 2\pi {\j}\int _{\zeta }\w\). This
is because he is using the branch of the logarithm with imaginary
part \(\in (0,2\pi) \). Also, we have used the homological version
of the Tate twist, which includes the factor \((2\pi\j)^{m}\).]
\end{example}

\section{\textbf{The map of complexes (second construction of 
\protect\( AJ\protect \))}}

We first describe the notation we shall use, which is a bit more 
involved than
that of the preceding section (there are also slight differences).

\subsection{Notation for Currents}

Let \( \X \) be a quasiprojective variety of complex dimension \( m 
\), \( Y\subset \X \)
an oriented analytic subset of real codimension \( k \), and \( \O 
\in \G (\O ^{\ell}_{\X }(\log D)) \)
where \( D\subset \X \) is any divisor. Associate to any given 
meromorphic
function \( f\in \C (\X ) \) the \( (2m-1) \)-chain \( 
T_{f}:=f^{-1}(\R ^{-}) \)
oriented so that \( \d T_{f}=(f)=|(f)_{0}|-|(f)_{\infty }|. \) Now 
define a current \( (\log f)\O \cdot \delta_{Y}\in F^{k}\DI^{\ell 
+k}(\X ) \) by

\begin{equation}
\int_{\X}(\log f) \O \cdot \delta_Y \v \w := \lim_{\e\to 0}\int_{Y\m 
\N_{\e}\left( D\cup T_f \right) }(\log f)\O \v \i^*_Y \w \\
\end{equation}
provided the limit exists for every \( C^{\infty } \)-form \( \w \in 
\G (\O ^{2m-\ell -k}_{\X ^{\infty }}) \)
compactly supported away from the boundary of \( \X \). (The \( F^{k} 
\)
means that all \( \w \in F^{m-k+1} \) are annihilated.)

\begin{rem*}
To re-iterate, in the r.h.s., ``\( \log f \)'' is always taken to 
have imaginary part \( \in (-\pi ,\pi ) \).
\end{rem*}
Recall that to any \( i \)-current \( \mathcal{K} \) is associated an \( (i+1) 
\)-current
\( \di [\mathcal{K}] \):
\[
\int _{\X }\di [\mathcal{K}]\v \w \: =(-1)^{i+1}\int _{\X 
}\mathcal{K}\v \di \w .\]
So for example \( \di [\log f]=\dlog f-2\pi \j\delta _{T_{f}} \) , 
and \( \di [\dlog f]=2\pi \j\delta _{(f)}. \)

\subsection{Higher Chow Groups}

We shall use the notation
\[
\s ^{n}\, :=\, (\P ^{1}_{\C }\m \{1\})^{n}\text {\, \, with\, 
coordinates\, }(z_{1},\ldots ,z_{n})\]
for affine \( n \)-space, with subsets
\[
\d \s ^{n}\, :=\, \bigcup _{i=1}^{n}\left\{ (z_{1},\ldots ,z_{n})\in 
\s ^{n}\, \left| \, z_{i}\in \{0,\infty \}\right. \right\} \, =\, 
\text {faces\, of\, }\s ^{n}\]
and \[
\d ^{k}\s ^{n}\, :=\, \bigcup _{i_{1}<\cdots <i_{k}}\left\{ 
(z_{1},\ldots ,z_{n})\in \s ^{n}\, \left| \, z_{i_{1}},\ldots 
,z_{i_{k}}\in \{0,\infty \}\right. \right\}
\]
\[\, =\, \text {codimension-}k\text {\, subfaces.}\]
Also let \[
\N _{\e }(\d \s ^{n})\, :=\]
\[\bigcup _{i=1}^{\infty }\left\{ (z_{1},\ldots ,z_{n})\in \s ^{n}\, 
\left| \, |z_{i}|<\e \text {\, or\, }|z_{i}|>\frac{1}{\e }\right. 
\right\} ,\text {\, \, \, and\, \, \, }\s ^{n}_{\e }:=\s ^{n}\m \N 
_{\e }(\d \s ^{n})\]

Let \( X \) be a complex projective variety of dimension \( m \), and 
define
subgroups of algebraic cycles on \( X\times \s ^{n} \)
\[
Z^{p}(X\times \s ^{n})\, \supseteq \, c^{p}(X,n)\, \supseteq \, 
d^{p}(X,n)\]
generated (resp.) by those subvarieties intersecting all subfaces \( 
X\times \d ^{k}\s ^{n} \)
properly, and (among those) by subvarieties pulled back from \( 
X\times \text {face} \)
by a coordinate projection. Set \( Z^{p}(X,n):=c^{p}(X,n)/d^{p}(X,n) 
\); writing
\( \rho _{i}^{0},\, \rho _{i}^{\infty } \) for the inclusions of the 
\( i^{\text {th}} \)
faces, define Bloch's differential
\[
\D :=\sum _{i=1}^{n}(-1)^{i-1}(\rho ^{\infty }_{i}{}^{*}-\rho 
^{0}_{i}{}^{*})\, :\, \, Z^{p}(X,n)\to Z^{p}(X,n-1).\]
Since \( \D \circ \D =0 \) this gives a complex, with \( CH^{p}(X,n) 
\) as
homology groups; only for our purposes cohomological indexing is 
better and
we shall write
\[
CH^{p}(X,n):=H^{-n}\left\{ Z^{p}(X,-\b )\right\} .\]

\subsection{Currents on \protect\( \s ^{n}\protect \)}

Set \[
\O ^{n}\, =\, \O (z_{1},\ldots ,z_{n})\, :=\, \dlog z_{1}\v \cdots \v 
\dlog z_{n}\, \in \, F^{n}\DI^{n}(\s ^{n}),
\]
holomorphic \(n\)-current.
\[
T^{n}\, :=\, T_{z_{1}}\cap \cdots \cap T_{z_{n}}\, \in \, \c _{n}(\s 
^{n})\text {\, topological\, }n\text {-chain}.\]
\[
R^{n}\, =\, R(z_{1},\ldots ,z_{n})\, :=\log z_{1}\dlog z_{2}\v \cdots 
\v \dlog z_{n}
\]
\[
+(\pm 2\pi \j)\log z_{2}\dlog z_{3}\v \cdots \v \dlog z_{n}\cdot 
\delta _{T_{z_{1}}}+\cdots \]
\[
\cdots +(\pm 2\pi \j)^{n-1}\log z_{n}\cdot \delta _{T_{z_{1}}\cap 
\cdots \cap T_{z_{n-1}}}\, \, \in \, \DI^{n-1}(\s ^{n}),\]
where ``\( \pm \)'' means \( (-1)^{n-1} \). For example, \( 
R^{1}=\log z \)
and \( R^{2}=\log z_{1}\dlog z_{2}-2\pi \j\log z_{2}\cdot \delta 
_{T_{z_{1}}}. \)
One may view these also as currents on \( X\times \s ^{n} \) by 
pullback.

From above one has for \( n=1 \) \( \di [R^{1}]=\O ^{1}-2\pi \j\cdot 
\delta _{T^{1}} \),
which generalizes to \( n>1 \) via\\
\begin{equation}
\di [R^n] = \O^n-(2\pi \j)^n \delta_{T^n} - (2\pi \j) 
\sum_{i=1}^n(-1)^i R(z_1,\ldots ,\hat{z_i},\ldots ,z_n) \cdot 
\delta_{(z_i)}.\\
\end{equation} Moreover,
\begin{equation}
\di[\O^n] = 2\pi \j \sum_{i=1}^n (-1)^i \O(z_1,\ldots 
,\hat{z_i},\ldots ,z_n) \cdot \delta_{(z_i)}
\end{equation}
and
\begin{equation}
\d T^n = \sum_{i=1}^n (-1)^{i} \left( \rho_i^0 {}_* T^{n-1} 
-\rho^{\infty}_i {}_* T^{n-1} \right). \\
\end{equation}

\subsection{Currents on \protect\( X\protect \)}

To produce these we must first specify a subcomplex of \( Z^{p}(X,\bullet) \)
consisting of elements in good position with respect to certain \( real \)
subsets of \( X \times \s^n \).  Let \( c^p_{\R}(X,n) \) consist of all \( \z \in 
c^p(X,n) \) intersecting \( X\times \left( T_{z_1}\cap\cdots\cap T_{z_j} \right) \)
and \( X \times \left\{ \left( T_{z_1}\cap \cdots \cap T_{z_j} \right) \cap \d^k\s^n \right\} \)
properly (for all \( 1\leq j \leq n \), \( 1\leq k <n \) ), and \( d^p_{\R}(X,n):=c^p_{\R}(X,n)\cap d^p(X,n) \).  
Then \( Z^p_{\R}(X,\bullet ):=c^p_{\R}(X,\bullet)/d^p_{\R}(X,\bullet) \) is a complex under \( \D \).  A moving 
technique based on unpublished notes of Bloch and worked out and extended by Levine in 
\cite[1.3.4]{Lev2} and \cite[Sect. 3.5.12]{Lev3} shows that this is quasi-isomorphic to the 
Bloch complex:

\subsection*{Moving by Translation Lemma:}

\( Z^p_{\R}(X,\bullet) \buildrel \simeq \over \to Z^p(X,\bullet). \)\\
\ \\
The proof consist of showing that any cycle is equivalent to one in general position 
after a generic complex affine translation in $\Delta^n$. 
Such cycles are obviously contained in $Z^p_{\R}(X,\bullet)$. 
Note that the arguments given in loc.cit. are using simplicial coordinates. But the quasiisomorphism between 
simplicial and cubical coordinate systems can be applied here. Alternatively one could restate 
the moving by translation lemma into cubical coordinates and prove it there.     
Let us offer an indication of what the "move" is.
Any \( \z \in Z^p(X,n) \) is defined over some \( k\subseteq \C \) finitely generated \( /\bar{\Q} \).  
Now consider \( \alpha_1, \ldots, \alpha_n \in \C^* \) such that \( trdeg \left( k(\alpha_1, 
\ldots, \alpha_n )/k \right) = n \), and let \( \tau=(\alpha_1,\ldots ,\alpha_n) \) act by 
multiplication on the coordinates \( (z_1,\ldots,z_n) \) of \( \z \) in \( \s^n \).

\begin{prop*}
Under these conditions, \( \tau \cdot \z \in Z^p_{\R}(X,n) \).
\end{prop*}

(Of course, \( \tau \cdot \z \) is no longer defined over \( k \).)

With some work, this proposition can be used to produce a map of complexes 
\( T: Z^p(X/k,\bullet)\to Z^p_{\R}(X,\bullet) \) together with a homotopy 
\( \mathcal{H}: Z^p(X/k,\bullet)\to Z^p(X,\bullet + 1) \) respecting \( Z^p_{\R} \) 
and subsets of \( X \), which satisfies \( T(\z)-\z = \D \mathcal{H}(\z) + \mathcal{H}(\D \z) \).  
The Lemma follows.

Associated to \( \z \in Z^{p}(X,n) \) one now produces \( \O _{\z 
}\in F^{p}\DI ^{2p-n}(X) \),
\( R_{\z }\in \DI ^{2p-n-1}(X) \) by the formulas
\[
\int _{X}\left\{ \begin{array}{c}
\O _{\z }\\
R_{\z }
\end{array}\right\} \v \w \, :=\, \lim _{\e \to 0}\sum _{j}n_{j}\int 
_{\z ^{j}_{\e }}\pi _{\s }^{\z _{j}}{}^{*}\left\{ \begin{array}{c}
\O ^{n}\\
R^{n}
\end{array}\right\} \v \pi _{X}^{\z _{j}}{}^{*}\w \]
where \( z=\sum n_{j}\z _{j} \) (\( \z _{j} \) irred.), \( \z 
_{j}^{\e }:=\z _{j}\cap (X\times \s _{\e }^{n}) \),
and \( \w \) is any \( C^{\infty } \)-form of the right 
degree.\footnote{%
That the limit on the r.h.s. always converges follows from an 
elementary analytic
argument, in which the proper intersection condition on each \( \z 
_{j} \)
plays a crucial role (e.g., see \cite{Ke} for the proof for \( \O 
_{\z } \)).
} These currents are zero if \( \z \in d^{p}(X,n). \)

\begin{rem*}
An appealing alternative form of the definition, e.g. for \( R_{\z } 
\), is
\( R_{\z }=\sum n_{j}\pi _{X}^{\z _{j}}{}_{*}\pi _{\s }^{\z 
_{j}}{}^{*}R^{n} \).
The push-forward \( \pi _{X}^{\z _{j}}{}_{*} \) should be regarded as 
involving
integration for those \( j \), for which \( \z _{j} \) has fibers of 
\( \dim \geq 1 \)
over \( X \). Note in particular that the numbers \( \text 
{codim}_{X}\left\{ \text {supp}\pi _{X}(\z _{j})\right\} \)
are \( not \) in general all the same.
\end{rem*}
Finally set \( T_{\z }:=\sum _{j}n_{j}\cdot \pi _{X}\left\{ \z 
_{j}\cap (X\times T^{n})\right\} \in \c _{2m-2p+n}(X). \)
The relations (5.2)-(5.4) give rise to formulas\\
\begin{equation}
\d T_{\z} = T_{\D \z}, \mspace{20mu} \di [\O_{\z}] = 2\pi \j \O_{\D 
\z}, \mspace{20mu} \di [R_{\z}] = \O_{\z} -(2\pi \j)^n 
\delta_{T_{\z}} -2\pi \j R_{\D \z}.
\end{equation}

\subsection{The map of complexes}

Define a complex of cochains for the Deligne homology of \( X \),
\[
\c _{\de }^{\b -2m}\left( X,\Z (p-m)\right) \, :=\, \Cone \left\{ 
\begin{array}{c} \c_{2m-\b }(X,\Z (p))\\
\bigoplus\\ F^{p}\DI ^{\b }(X)\end{array}\to \DI ^{\b }(X)
\right\} [-1](-m)\]
\[
=\left\{ \c _{2m-\b }(X,\Z (p))\oplus F^{p}\DI ^{\b }(X)\oplus \DI 
^{\b -1}(X)\right\} (-m)\]
with differential \( D \) taking \( (a,b,c)\mapsto (-\d a,-\di 
[b],\di [c]-b+\delta _{a}). \)
Then according to the formulas (5.5), sending \[
\z \mapsto \frac{(-2\pi \j)^{p-n}}{(2\pi \j)^{m}}\left( (2\pi 
\j)^{n}T_{\z },\O _{\z },R_{\z }\right) =:\r _{X}(\z )\]
produces a map of complexes
\[
\r _{X}:\, Z^{p}_{\R}(X,-\b )\to \c _{\de }^{2p-2m+\b }(X,\Z (p-m));\]
that is, \( D\r _{X}(\z )=\r _{X}(\D \z ). \) (Note: \( \O _{\z }=0 
\) if
\( p>m \) or \( p<n \).) According to the Moving Lemma we may replace \( Z^p_{\R}(X,-\bullet) \) by \( Z^p(X,-\bullet) \) with the \( caveat \) that the map is in the derived category.  This induces the desired map
\[
AJ:\, CH^{p}(X,n)\to H^{\de }_{2m-2p+n}(X,\Z (m-p))\begin{array}{c}
_{\cong }\\
\longleftarrow \\
^{\text {P.D.}}
\end{array}H_{\de }^{2p-n}(X,\Z (p)).\]
If \( \D \z =0 \) then \( \z \) represents a class \( [\z ]\in 
CH^{p}(X,n) \),
and we write \( AJ[\z ] \) or \( \Phi_{p,n}(\z) \) for the class \( 
[\r
_{X}(\z )] \).

\begin{rem*}
(i) \( AJ \) is in fact a ring homomorphism; that is, if \( [\W ]\in 
CH^{p}(X,\ell ) \)
and \( [\Y ]\in CH^{q}(X,n) \) then \( [\W \times \Y ]\in 
CH^{p+q}(X,\ell +n) \)
and \[
[\r _{X}(\W )]\cup [\r _{X}(\Y )]=[\r _{X}(\W \times \Y )]\]
under the cup-product in Deligne (co)homology. The class on the 
l.h.s. is given
(modulo factors of \( 2\pi \j \)) by
\[
\left( (2\pi \j)^{\ell +n}T_{\W }\cap T_{\Y },\O _{\W }\v \O _{\Y 
},(-1)^{\ell }(2\pi \j)^{\ell }\delta _{T_{\W }}\cdot R_{\Y }+R_{\W 
}\v \O _{\Y }\right) ;\]
that this equals \( \left( (2\pi \j)^{\ell +n}T_{\W \times \Y },\O 
_{\W \times \Y },R_{\W \times \Y }\right) \)
is implied by the formula
\[
R(w_{1},\ldots ,w_{\ell };y_{1},\ldots ,y_{n})
\]
\[
=(-1)^{\ell }(2\pi i)^{\ell }\delta _{T(w_{1},\ldots ,w_{\ell 
})}\cdot R(y_{1},\ldots ,y_{n})\, +\, R(w_{1},\ldots ,w_{\ell })\v \O 
(y_{1},\ldots ,y_{n})\]
on \( \s ^{\ell +n}
=\s ^{\ell }\times \s ^{n} \) (with coordinates \( w_{1},\ldots 
,w_{\ell };y_{1},\ldots ,y_{n} \)).

(ii) The projection of this \( AJ \) map to the real Deligne 
cohomology, i.e.
the composition
\[
CH^{p}(X,n)\begin{array}[t]{c}
\longrightarrow \\
^{AJ}
\end{array}H_{\de }^{2p-n}(X,\Z (p))\begin{array}[t]{c}
\longrightarrow \\
^{\pi _{\R }}
\end{array}H_{\de }^{2p-n}(X,\R (n))\]
agrees exactly with the regulator map defined by Goncharov in 
\cite{Go}\footnote{%
The \( AJ \) map defined there, on the other hand, was not correct.
} (see \cite{Ke}).
\end{rem*}

\subsection{Passage to ordinary cohomology}

Now let \( n\geq 1 \) and \( \z \) be a higher Chow cycle: \( \D \z 
=0 \).
Then \( \di [\O _{\z }]=0 \), \( \d T_{\z }=0 \), and \( \di [R_{\z 
}]=\O _{\z }-(2\pi \j)^{n}T_{\z }\, \, \implies \, \, [\O _{\z 
}]=(2\pi \j)^{n}T_{\z } \) in \( H^{2p-n}(X,\C) \) .
Multiplying by  \( (2\pi \j)^{p-n} \), we get a class in \[
 (*) \, \, \, \, \, \, \, \, \, F^{p}H^{2p-n}(X,\C )\cap H^{2p-n}(X,\Z
(p)).\]
This is the closest thing we get to a ``fundamental class'' for \( \z 
\) ; when \( T_{\z } \sim 0 \) we say \( [\z ] \in CH^p_{hom}(X,n) \).

For \( X \) projective (and \( n\geq 1 \) ), \( (*) \) is a torsion
group.  Consequently, the class \( [ \Omega_{\z} ] \in F^p H^{2p-n}
(X,\C) \) is zero while \( [T_{\z} ] \in H^{2p-n}(X,\Z) \) is at worst
torsion.

To proceed further we must have \( T_{\z} \sim 0 \).  So in general,
we must either (a) pass to rational coefficients (to render \( (*) \) 
zero) or (b) assume the slight restriction \( [\z ] \in
CH^p_{hom}(X,n) \).  While we have chosen (b), we emphasize that what 
follows (for the remainder of the paper) works essentially verbatim with
\( \Z \) replaced everywhere by \( \Q \) (instead of this assumption).
Moreover, if \( X \) is such that \( H^{2p-n}(X,\Z) \) has no torsion, \(
CH^p(X,n) = CH^p_{hom}(X,n) \) and no such choice is neccessary.  (Trivial
example:  \( X = \) pt.)
 
Assuming, then, that \(\z\) (i.e., \( T_{\z} \) ) is nullhomologous,
there exist ``primitives'' \( \Xi \in F^{p}\DI ^{2p-n-1}(X) \),
\( \zeta \in \c _{2m-2p+n+1}(X,\Z) \) such that \( \di [\Xi ]=\O _{\z } 
\), \( (-1)^{n}d[\delta_{\zeta}] = \delta_{\d \zeta} = T_{\z } \) (or,
strictly speaking, \( \delta_{T_{\z}} \) ).  \( \zeta \) is called a
``membrane''.

We may now modify \( \r_{X}(\z ) \) by a coboundary, to get
\[
\left( (2\pi \j)^{n}T_{\z },\O _{\z },R_{\z }\right) +D\left( (-2\pi 
\j)^{n}\zeta ,\Xi ,0\right) =\left( 0,0,R_{\z }-\Xi +(-2\pi 
\j)^{n}\delta _{\zeta }\, \, =:R_{\z }''\right) .\]
Now \( \Xi \) and \( (-2\pi \j)^{p}\zeta \) are ambiguous by \( 
F^{p}H^{2p-n-1}(X,\C ) \)
and \( H^{2p-n-1}(X,\Z (p)) \), respectively, and so we have a 
well-defined
class
\[
(-2\pi \j)^{p-n}[R_{\z }'']\in \frac{H^{2p-n-1}(X,\C 
)}{F^{p}H^{2p-n-1}(X,\C )+H^{2p-n-1}(X,\Z (p))}\]
reflecting the isomorphism of the latter group with \( H_{\de 
}^{2p-n}(X,\Z (p)) \)
for \( n\geq 1 \) and \( X \) projective. {[}Note: \( (2\pi \j)^{m} 
\) has
already disappeared in the P.D. \( \cong \).{]} Since this quotient 
is equivalent
to \[ \left. \left\{ F^{m-p+1}H^{2m-2p+n+1}(X,\C )\right\} ^{\vee 
}\right/ \text {im}\left\{ H_{2m-2p+n+1}(X,\Z (p))\right\} , \]
this class (and thus \( AJ[\z ] \)) is computed by the functional
\[
\frac{1}{(-2\pi \j)^{n-p}}\int _{X}R_{\z }'\v (\cdot )\mspace 
{50mu}\text {modulo\, periods\, }(2\pi \j)^{p}\int _{\gamma }(\cdot 
)\]
where \( R_{\z }':=R_{\z }+(2\pi \j)^{n}\delta _{\zeta } \) and we 
may drop the \( \Xi \)-term, again thanks to Hodge-type considerations
(see remark 2.3).

\begin{example*}
One easily recovers Levine's formula for \( CH^{p}(X,1) \) from this 
approach,
writing \( R^{1}=\log z \) and \( R_{\z }=\sum _{\alpha }\pi 
_{X}^{(\z _{\alpha },f_{\alpha })}{}_{*}\pi _{\s }^{(\z _{\alpha 
},f_{\alpha })}{}^{*}(\log z)\, =\, \sum _{\alpha }\log f_{\alpha 
}\cdot \delta _{\z _{\alpha }}. \)
\end{example*}

\subsection{A further simplification for \protect\( n\geq p\protect 
\) or \protect\( p>m\protect \).}

For \( p \) in this range we may clean up the above functional 
considerably;
since then \( F^{p}H^{2p-n-1}(X,\C )=0 \) and \( F^{m-p+1}\O 
_{X^{\infty }}^{2m-2p+n+1}(X)=F^{0}\O ^{2m-2p+n+1}_{X^{\infty }}(X), 
\)
we may as well evaluate it on integral classes --- that is, 
Poincar\'e duals
of topological \( (2p-n-1) \)-cycles \( \xi \). The ``periods'' are 
then
all in \( \Z (p) \), as is the contribution from the term \( (2\pi 
i)^{n}\delta _{\zeta } \)
in \( R_{\z }' \). One may therefore regard the functional
\[
\frac{1}{(-2\pi \j)^{n-p}}\int _{(\cdot )}R_{\z }\, \, \, \in \, \, 
\, \text {Hom}\left( H_{2p-n-1}(X,\Z ),\left. \C \right/ \Z 
(p)\right) \]
as representing \( AJ[\z ]\in H^{2p-n-1}(X,\left. \C \right/ \Z (p)) 
\).

\begin{example*}
In the very simple (but interesting) series of examples \( 
CH^{p}(\text {pt.},2p-1), \)
\( R_{\z } \) is simply the number \( \int _{\z }R^{2p-1} \), and the 
resulting
classes
\[
AJ[\z ]=\frac{1}{(-2\pi \j)^{p-1}}\int _{\z }R^{2p-1}\, \, \in \, \, 
\left. \C \right/ \Z (p)\cong H^{1}_{\de }(\text {pt.},\Z (p))\]
are related to Goncharov's Chow \( p \)-logarithm \cite{Go}. Here we 
just
want to point out (for \( p=2 \) ) how to construct  
classes in \( CH^2(\text{pt.},3) \) with torsion and nontorsion
\( AJ \) images.

Consider \( a,\, b\, \in \C ^{*}\m \{1\} \)
and introduce \[
V(a)=\left\{ (1-\frac{a}{t},\, 1-t,\, t)\, \left| \, t\in \P 
^{1}\right. \right\} \cap \s ^{3},\]
\[
W(b)=\left\{ (1-\frac{b}{t},\, t,\, 1-t)\, \left| \, t\in \P 
^{1}\right. \right\} \cap \s ^{3}.\]
Then one easily shows
\[
\partial V(a)=(1-a,a),\, \, \, \, \partial W(b)=(b,1-b).\]
Hence \( \xi _{a}:=V(a)-W(1-a) \) defines a class in \( 
CH^{2}(\text {pt.},3). \)

The value \( \Phi _{2,3}(\xi _{a})\in \C /\Z (2) \) is not hard to 
compute.
Indeed one finds that
\[
\Phi _{2,3}(\xi _{a})=\text {Li}_{2}(a)+\text {Li}_{2}(1-a)+\log 
a\log (1-a),\]
where \( \text {Li}_{2} \) is the dilogarithm, and \( \log \) is the 
principal
branch. By Beilinson's rigidity, this value does not depend on \( a 
\). Hence
\[
\Phi _{2,3}(\xi _{a})=\lim _{a\to 0}\Phi _{2,3}(\xi _{a})=\sum 
_{n=1}^{\infty }\frac{1}{n^{2}}=\frac{\pi ^{2}}{6}\in \C /\Z (2)\]
is a torsion class.

Next, let \( D_2 = \) the Bloch-Wigner function, \( \mathcal{B}_2(\C) 
= \)
the Bloch group, and \( st: \mathcal{B}_2(\C) \to \C^* \bigwedge_{\z} 
\C^*
\) the standard map \( \left\{ a \right\} _2 \mapsto (1-a) \wedge a\).
Finally set \(
\rho (a) = Alt_3(V(a)) \). The following is proved in \cite{Ke} \( \S
3.1.2 \):
\begin{prop*}
Given any element \( \sum m_j \left\{ a_j \right\} _2 \in 
\text{ker}(st)
\subseteq
\mathcal{B}_2(\C) \), \( \sum m_j \rho(a_j) \in Z^2(\text{pt.},3) \) 
may
be completed to a higher Chow cycle \( \z \) by adding "decomposable"
elements \( \in Z^1(\text{pt.},2)\bigwedge Z^1(\text{pt.},1) \). The
composition \( \pi_{\R} \circ \Phi_{2,3} \) on \( \z \) is then 
computed
by \( \Im (R_{\z}) = \sum m_j D_2(a_j) \in \R \). \end{prop*}
So if \( \sum m_j D_2 (a_j) \neq 0 \) (there are many examples), \( 
AJ(\z)
\in \C/\Z(2) \) is nontorsion.
\end{example*}

\begin{example*}
Let \( X \) be a compact Riemann surface. We recover the formula for 
the real
regulator
\[
r_{2,2}:\, CH^{2}(X,2)\to H_{\de }^{2}(X,\R (2))\]
in \cite{Ra} by composing \( AJ \) with \( \pi _{\R } \) (which takes 
the
imaginary part in this case).

The irreducible components of \( \z \in Z^{2}(X,2) \) are of two 
types: (a)
curves \( \subset \s ^{2} \) over isolated points of \( X \); (b) 
graphs (over
\( X \)) of pairs of meromorphic functions \( f,\, g\, \in \C (X) \). 
Writing
\[
\Gamma _{f,g}:=\left\{ \left. (x,f(x),g(x))\in X\times (\P 
^{1})^{2}\, \, \right| \, x\in X\right\} \cap (X\times \s ^{2}),\]
one has \[
\z =\sum m_{k}\cdot \{x_{k}\}\times C_{k}\, +\, \sum n_{j}\Gamma 
_{f_{j},,g_{j}}\, =\, \z _{(a)}+\z _{(b)}.\]
Now assume \( \partial _{\mathcal{B}}\z =0 \), and compute \( AJ[\z ] 
\):
since \( \int _{C_{j}}R^{2}\v \pi _{X}^{*}\omega =0 \) (for any \( 1 
\)-form
on \( X \)), \( R_{\z _{(a)}}=0 \). Thus
\[
R_{\z }=R_{\z _{(b)}}=\sum n_{j}R_{f_{j},g_{j}}=\sum n_{j}(\log 
f_{j}\dlog g_{j}-2\pi \j\log g_{j}\delta _{T_{f_{j}}})\, \, ;\]
integrating the latter current over loops
\( \gamma \in \Omega X \) gives
\[
AJ[\z ]\in \text {Hom}\left( H_{1}(X,\Z ),\C /\Z (2)\right) \cong 
H^{1}_{\de }(X,\Z (2)).\]
The technical point here is: prior to integrating one should move the 
loop (in its homology
class) to avoid all points \( \{x_{k}\} \) and \( \bigcup 
|(f_{j})|\cup |(g_{j})| \).
It can be shown directly (and is also clear from the setup here) that 
integrals
over loops around any of these points are trivial (\( \in \Z (2) \)), 
and so
the \( \int \) is independent of the ``move'' in question.
Note that the imaginary part of \( \log f\dlog g-2\pi \j\log g\, 
\delta _{T_{f}} \)
is \( \log |f|\text {darg}g+\arg f\dlog |g|-2\pi \log |g|\delta 
_{T_{f}} \).
Adding to this \( \text {d}[-\arg f\log |g|]\, =\, -\log |g|\text 
{darg}f-\arg f\dlog |g|+2\pi \log |g|\delta _{T_{f}} \)
gives the cohomologous current\footnote{%
It is well-known that Beilinson refers to this as ``Mama's formula''.}
\[
\log |f|\text {darg}g-\log |g|\text {darg}f.\]
Integrating this over loops gives \[
r_{2,2}[\z ]\in \text {Hom}(H_{1}(X,\Z ),\R )\cong H^{2}_{\de }(X,\R 
(2)).\]
An alternate formula sends
\[
\gamma \mapsto \Im \left( \int _{\gamma }\log f\frac{dg}{g}-\log 
|g(p)|\int _{\gamma }\frac{df}{f}\right) \]
where \( \gamma \) is based at \( p \) and \( \log f \) is \( not \) 
the
principal branch (but rather is continued along \( \gamma \), 
starting from
\( p \)).
\end{example*}

\subsection{Geometric interpretation}

Now we ask, to what extent is \( AJ \) (for \( n\geq 1 \)) ``like'' 
an Abel-Jacobi
map in the classical sense (for \( n=0 \)), of integrating forms over 
a chain
\( \G \) with \( \d \G =\z \)? To answer this question we first 
extend the
classical approach to the subgroup \( Z^{p}(X\times \s ^{n},X\times 
\d \s ^{n})\subseteq c^{p}(X,n) \)
of \( relative \) (algebraic) cycles, consisting of those \( \z \) 
for which \( \z \cdot (X\times \d \s ^{n})=0 \),
i.e. the intersections with \( each \) \( face \) cancel (counted 
with multiplicity). Since every class in \( CH^{p}(X,n) \)
is represented by such a cycle, it makes sense to ask whether the 
resulting \( relative \) \( AJ \) \( map \) coincides with the \( AJ 
\) constructed above.

The relevant details in the discussion that follows can be found in 
\cite{Ke}.
Let \( \z \in Z^{p}(X\times \s ^{n},X\times \d \s ^{n}) \) have 
{[}complex{]}
dimension \( d=m+n-p \).  Assume as above that \( [\z ] \in
CH^p_{hom}(X,n) \), and for the time being that \( \z \) is also
in \( Z^p_{\R}(X,n) \).  We replace \( \z \) by a ``limit'' of 
topological
cycles via a kind of excision. Namely, writing \( \I ^{n}:=\bigcup 
_{i=1}^{n}\left\{ \left. (z_{1},\ldots ,z_{n})\in (\C ^{*})^{n}\, 
\right| \, z_{i}=1\right\} \)
and letting \( \z ^{0}_{\e } \) represent the analytic closure of \( 
\z \cap (X\times \s ^{n}_{\e }) \)
on \( X\times (\C ^{*})^{n} \), for each \( \e >0 \) (suff. small) 
there are
topological relative cycles
\[
\z ^{0}_{\e }+\W _{\e }=\z _{\e }\in Z_{2d}^{\text {top}}\left( 
X\times (\C ^{*})^{n},X\times \I ^{n}\right) \]
where \( \lim _{\e \to 0}\int _{\W _{\e }}\pi _{\s }{}^{*}\left\{ 
\begin{array}{c}
\O ^{n}\\
R^{n}
\end{array}\right\} \v \pi _{X}{}^{*}\alpha \, \, \, =\, 0 \) (\( 
\forall \) \( C^{\infty } \)-forms \( \alpha \) on \( X \)). Since
Lefschetz duality guarantees a perfect pairing between
\[
H_{2d}\left( X\times (\C ^{*})^{n},X\times \I ^{n}\right) \cong 
H_{2d-n}\left( X\right) \otimes \left\langle (S^{1})^{n}\right\rangle 
\]
and
\[
H_{2p}\left( X\times \s ^{n},X\times \d \s ^{n}\right) \cong 
H_{2p-n}\left( X\right) \otimes \left\langle T^{n}\right\rangle ,\]
\( \z _{\e } \) is homologous (mod \( X\times \I ^{n} \)) to \( T_{\z 
_{\e }}\times (S^{1})^{n} \),
where \( T_{\z _{\e }}=\pi _{X}\left( \z _{\e }\cap (X\times 
T^{n})\right) \).
So there exists a topological \( (2d+1) \)-chain \( \G ^{0}_{\e } \) 
on \( X\times (\C ^{*})^{n} \)
such that
\[
\d \G _{\e }^{0}-\z _{\e }+T_{\z _{\e }}\times (S^{1})^{n}\, \, 
\subseteq \, \, X\times \I ^{n}.\]
The relative cycle \( \z \) is also a higher Chow cycle; therefore \( 
T_{\z _{\e }}[\sim T_{\z }]\sim 0 \)
is the boundary of a membrane \( \zeta _{\e } \) on \( X \). (This is 
equivalent
to \( [\z ]\in CH^{p}_{\text {hom}}(X,n) \), which is
what we are assuming.)

Set \( \G _{\e }=\G _{\e }^{0}+\zeta _{\e }\times (S^{1})^{n}. \) 
Since this
has boundary \( \z _{\e } \) (mod \( X\times \I ^{n} \)), it now 
makes sense
to define the relative Abel-Jacobi of \( \z \) as a functional (mod 
periods)
on
\[
F^{d+1}H^{2d+1}\left( X\times (\C ^{*})^{n},X\times \I ^{n};\C 
\right) \cong F^{d-n+1}H^{2d-n+1}(X,\C )\otimes \left\langle 
\frac{1}{(2\pi \j)^{n-p}}\O ^{n}\right\rangle \]
induced by the formula
\[
AJ_{\text {rel}}(\z )\w \, :=\, \lim _{\e \to 0}\frac{1}{(2\pi 
\j)^{n-p}}\int _{\G _{\e }} \pi _{\s }{}^{*}\O^{n} \v\pi _{X}{}^{*}\w
\]
where as test forms we use all \( \di \)-closed \( \w \in F^{m-p+1}\O 
^{2m-2p+n+1}_{\X ^{\infty }}(X). \)
It can be shown that the resulting map
\[
AJ_{\text {rel}}:\, Z^{p}(X\times \s ^{n},X\times \d \s ^{n})\to 
H_{\de }^{2p-n}(X,\Z (p))
\]
\[
\hskip1.5in \cong \frac{H^{2p-n-1}(X,\C )}{F^{p}H^{2p-n-1}(X,\C 
)+H^{2p-n-1}(X,\Z (p))}\]
respects relative rational equivalence. (Intuitively speaking, \( 
AJ_{\text {rel}} \) should go to some ``\( H_{\de }^{2p}\left( 
X\times \s ^{n},X\times \d \s ^{n};\Z (p)\right) \)'';
one can easily justify defining this to \( be \) \( H_{\de 
}^{2p-n}(X,\Z (p)) \).)

To see this is the same as \( AJ \) of \( \z \) (considered instead 
as a
higher Chow cycle), we show
\begin{prop}\label{currentequate}
\[
\lim _{\e \to 0}\int _{\G _{\e }}\O^{n}\v\w\, =\, (-1)^{n}\biggl(\int 
_{\z }
R^{n}\v\w\, +\, (-2\pi \j)^{n}\int _{\zeta }\w\biggr) \]
for \( \w \) \( \di \)-closed. \end{prop}
\begin{proof}
By Stokes' theorem for currents and equation \( (5.2) \), \[
\int _{\z _{\e }(=\d \G _{\e })}R^{n}\v\w\, =\, (-1)^{n}\int _{\G 
_{\e }}\di [R^{n}\v\w]\, =\, (-1)^{n}\int _{\G _{\e }}\di 
[R^{n}]\v\w\, =\, \]
\[
(-1)^{n}\int _{\G _{\e }}\O ^{n}\v\w \, -\, (-2\pi \j)^{n}\int _{\pi 
_{X}[\G _{\e }\cap (X\times T^{n})]}\w \, \, \, -\, 0.\]
The residue term in equation \( (5.2) \) makes no contribution 
because \( \G _{\e }\cap (X\times \d \s ^{n})=\emptyset \).
Now one can construct \( \G _{\e }^{0} \) explicitly (as in 
\cite{Ke}) so
that \[ \dim _{\R }[\pi _{X}\{\G ^{0}_{\e }\cap (X\times 
T^{n})\}]<2d-n+1 ;\]
for our purposes then \( \pi _{X}[\G _{\e }\cap (X\times 
T^{n})]=\zeta _{\e } \),
and
\[
\int _{\G _{\e }}\w \v \O ^{n}\, =\, (-1)^{n}
\biggl(\int _{\z _{\e }}\w \v R^{n}\, +\, (-2\pi \j)^{n}\int _{\zeta 
_{\e }}\w\biggr) .\]
Taking limits then gives the result.
\end{proof}

\begin{rem*}
(i) In light of the coincidence of these two maps, we may view the 
previous ``simplification
for \( n\geq p \) or \( p>m \)'' as asserting that the relative \( AJ 
\)
map on \( (X\times \s ^{n},X\times \d \s ^{n}) \) may be ``pushed 
down to
\( X \)'' where it consists merely of computing periods of the 
current \( R_{\z } \).
This is essentially thanks to the fact that \( H^{n}((\C ^{*})^{n},\I 
^{n})=F^{n}H^{n}((\C ^{*})^{n},\I ^{n}) \),
and \( H^{*}((\C ^{*})^{n},\I ^{n})=0 \) for \( *\neq n \).

(ii) We can modify this approach in case \( \z \notin Z^p_{\R}(X,n) \).
Referring to the Proposition in \( \S 5.4 \), we remark that \( \z \) can also be brought into good position with respect to the \( X\times \left( T_{z_1}\cap \cdots \cap T_{z_i} \right) \) etc. by perturbing \( these \) \( real \) \( chains \) (rather than \( \z \) ) by \( \tau \).

To deal with such \( \z \), therefore, one merely repeats all of \( \S 5.8 \) with the "perturbation" \( T_{z_i} \mapsto T_{z_i}' := T_{z_i/\alpha_i} \) (which accordingly leads to a different chain \( \Gamma_{\epsilon}' \) ), and also with branches of \( log(z_i) \) in \( R^n \) replaced by branches with cuts at \( T_{z_i}' \).  Proposition \( 5.1 \) then holds exactly.  There is no need to take limits as the perturbations approach \( 1 \), since \( \Gamma_{\epsilon}' \) amounts simply to a different choice of bounding membrane for computing the relative AJ map.

Such a procedure has been carried out in a concrete computation in \( \S 3.2 \) of \cite{Ke}.
\end{rem*}

\subsection{Quasiprojective case}

Now let \( V\subset X \) be an arbitrary divisor; we show how to 
define \[
AJ_{X\m V}:\, CH^{p}(X\m V,n)\to H_{\de }^{2p-n}(X\m V,\Z (p)).\]
Let \( (\bar{X},\bar{V}) \) be a pair where \( \bar{V} \) is a n.c. 
divisor
and \( \bar{X}\m \bar{V}=X\m V \). By Bloch's moving lemma, 
restriction induces
a quasiisomorphism of complexes
\[
\left. Z^{p}(\bar{X},-\b )\right/ Z^{p-1}(\bar{V},-\b 
)\begin{array}[t]{c}
\longrightarrow \\
^{\simeq }
\end{array}Z^{p}(X\m V,-\b ).\]
The homotopy \( \mathcal{H} \) described in \( \S 5.4 \) leads to a proof that
\[
\left. Z^p_{\R}(\bar{X},-\b)\right/ Z^{p-1}_{\R}(\bar{V},-\b) \buildrel \simeq \over \to 
\left. Z^p(\bar{X},-\b)\right/ Z^{p-1}(\bar{V},-\b). \]
Finally, the triple \( \left( (2\pi \j)^{n}T_{\z },\O _{\z },R_{\z }\right) \) 
once
again yields a map of complexes
\[
AJ_{X\m V}:\, \left. Z^{p}_{\R}(\bar{X},-\b )\right/ Z^{p-1}_{\R}(\bar{V},-\b 
)\longrightarrow \]
\[
\Cone \left\{ \begin{array}{c}\c _{2m-2p-\b }(\bar{X},\bar{V};\Z 
(p))\\
\bigoplus\\ \G \left( F^{p}\DI _{\bar{X}}^{2p+\b }(\log 
\bar{V})\right)\end{array} \to \G \left( \DI _{\bar{X}}^{2p+\b }(\log 
\bar{V})\right) \right\} [-1](-d),\]
which induces \( AJ_{X\m V}. \)

\begin{rem*}
(i) We need the n.c. condition in order that the latter complex 
actually compute
\( H^{2p-n}_{\de }(X\m V,\Z (p)) \).

(ii) We could not use \( Z^{p}(\bar{X}\m \bar{V},-\b ) \) here 
because \( \O _{\z } \)
and \( R_{\z } \) have (in general) worse than log poles along \( 
\bar{V} \),
for \( \z \in Z^{p}(\bar{X}\m \bar{V},n). \)

(iii) the simplifications that occurred in the projective case for \( 
n \geq 1 \) (and rational coefficients),
require \( p>m \) or \( n>p \) here (as \( F^{p}H^{2p-n}(X\m V,\C 
)\cap H^{2p-n}(X\m V,\Q (p)) \)
must vanish).

(iv) In case \( n>p, \) or \( n=p>m, \) the \( \G \left( F^{p}\DI 
_{\bar{X}}^{2p+\b }(\log \bar{V})\right) \)
terms do not enter. The situation simplifies and we may work with \( 
R_{\z } \)
(resp. \( T_{\z } \)) directly on \( X\m V \) (resp \( (X,V) \)). 
Taking
the limit we get \( AJ \) maps over the generic point
\[
AJ_{\eta _{X}}:\, CH^{p}(\C (X),n)\to H_{\de }^{2p-n}(\eta _{X},\Z 
(p)).\]

\end{rem*}
\begin{example*}
For \( n=p>m \), \( CH^{n}(\C (X),n)\cong K_{n}^{M}(\C (X)) \) by a 
result
of Totaro\footnote{%
also due to Nesterenko and Suslin
} \cite{To}. The resulting map
\[
AJ_{\eta _{X}}:\, K_{n}^{M}(\C (X))\to H^{n}_{\de }(\eta _{X},\Q 
(n))\cong H^{n-1}(\eta _{X},\C /\Q (n))\]
is called the Milnor regulator and is studied extensively in 
\cite{Ke}. As
in the projective case, \( AJ_{\eta _{X}}\{f_{1},\ldots ,f_{n}\} \) 
may be
computed as a functional on topological cycles, namely
\[
\int _{(\cdot )}R(f_{1},\ldots ,f_{n})\, \, \in \, {\rm Hom}\left( 
H_{n-1}(\eta _{X},\Z ),\, \C /\Q (n)\right) .\]

\end{example*}

\section{\textbf{How a cycle gives rise to an extension of motives \\ 
(third construction
of \protect\( AJ\protect \))}}

Let \( \z \) be a cycle in \( CH_{\text {hom}}^{p}(X,n) \). By 
normalization
of chain complexes (\cite[page 512]{E-Z}), we may assume that all 
individual intersections
\( \z \cap \{z_{i}=0,\infty\} \) are zero. Let \( U:=\s _{X}^{n}\m |\z | \) 
and \( \partial U:=U\cap \partial \s _{X}^{n}. \)
One has an exact sequence (in any reasonable theory satisfying weak 
purity and
the homotopy axiom)
\[
H^{2p-2}(U)\to H^{2p-2}(\partial U)\to H^{2p-1}(U,\partial U)\to 
H^{2p-1}(U)\to H^{2p-1}(\partial U).\]
By weak purity \( H^{i}_{|\z |}(\s _{X}^{n})=0 \) for \( i<2p \). 
Also by
the homotopy axiom, \( H^{i}(\s ^{n}_{X})=H^{i}(X) \) for \( i\geq 0 
\). Therefore
\( H^{2p-2}(U)=H^{2p-2}(X) \) (naturally) and furthermore one has for 
all \( i\geq 0 \)
\[
H^{i}(\partial \s _{X}^{n})=H^{i}(X)\oplus H^{i-n+1}(X),\]
i.e. \( \partial \s ^{n}_{X} \) is like a real \( (n-1) \)-sphere. 
Moreover
\( H^{i}(\partial U)=H^{i}(\partial \s ^{n}_{X}) \) for \( i=2p-2 \), 
again
by weak purity for the faces. Using all this, the long exact sequence 
now becomes
\[
0\to H^{2p-n-1}(X)\to H^{2p-1}(U,\partial U)\to \ker \left\{ 
H^{2p-1}(U)\to H^{2p-1}(\partial U)\right\} \to 0.\]
But \[
\ker \left\{ H^{2p-1}(U)\to H^{2p-1}(\partial U)\right\} \subseteq 
\ker \left\{ H_{|\z |}^{2p}(\s _{X}^{n})^{\circ }\begin{array}[b]{c}
_{\beta }\\
\to \end{array}H_{|\partial \z |}^{2p}(\partial \s _{X}^{n})^{\circ 
}\right\} .\]
The symbol \( ^{\circ } \) stands for the kernel of the map 
forgetting supports.
This implies that we have a long exact sequence

\begin{equation}
0 \to H^{2p-n-1}(X) \to H^{2p-1}(U,\partial U) \to \ker (\beta) \to 
H^{2p-n}(X). \\
\end{equation}
To see this, one simply applies the serpent lemma to the diagram of 
exact sequences:

\[ \begin{array}{ccccccccc}
& & 0 & & H^{2p-2}(\partial U) \\
& & \downarrow & & \downarrow\\
& & H^{2p-1}(X) & & H^{2p-1}(U,\partial U) \\
& & || & & \downarrow & \searrow \\
0 & \to & H^{2p-1}(\s^n_X) & \to & H^{2p-1}(U) & \to & 
H^{2p}_{|\z|}(\s^n_X)^{\circ} & \to & 0 \\
& & \downarrow & & \downarrow & & \downarrow \\
0 & \to & H^{2p-1}(\partial \s^n_X) & \to & H^{2p-1}(\partial U) & 
\to & H^{2p}_{|\partial \z|}(\partial \s^n_X)^{\circ} & \to & 0 \\
& & \downarrow \\
& & H^{2p-n} (X) \\
& & \downarrow \\ & & 0 \\
\end{array} \]
This yields an extension via pullback

\begin{equation}
0 \to H^{2p-n-1}(X,\Z(p)) \to \mathbf{E} \to \Z(0) \to 0 , \\
\end{equation}
where the motive \( \Z (0) \) is generated by the algebraic cycle \( 
\{\z \} \)
in \( \ker (\beta ) \). If we specialize to singular cohomology, we 
obtain
an extension of mixed Hodge structures:

\begin{equation}
\mathbf{E} \in \text{Ext}^1_{\text{MHS}}(\Z(0), H^{2p-n-1}(X,\Z(p))) 
\end{equation}
\[
\hskip1.6in
\simeq \frac{H^{2p-n-1}(X,\C)}{F^p H^{2p-n-1}(X,\C) + 
H^{2p-n-1}(X,\Z(p))} =: J^{p,n}(X) , \\
\]
as desired.

Note that \( H_{|\z |}^{2p}(\s ^{n}_{X}) \) is generated by the 
components
of \( |\z | \), and therefore \( F^{p}H_{|\z |}^{2p}(\s 
^{n}_{X})=H_{|\z |}^{2p}(\s ^{n}_{X}). \)
Now put
\[
V=\ker \left\{ \ker (\beta )\to H^{2p-n}(X)\right\} .\]
Then \( F^{p}V=V \). By applying the serpent lemma to

\[
\begin{array}{ccccccccc}
0 & \to & F^p H^{2p-n-1}(X) & \to & F^p H^{2p-1}(U,\partial U) & \to 
& F^p V & \to & 0 \\
&\\
& & \downarrow & & \downarrow & & || \\
&\\
0 & \to & H^{2p-n-1}(X) & \to & H^{2p-1} (U,\partial U) & \to & V & 
\to & 0 \\
\end{array}
\]
we deduce that
\[
\frac{H^{2p-n-1}(X,\C )}{F^{p}H^{2p-n-1}(X,\C )}\simeq 
\frac{H^{2p-1}(U,\partial U,\C )}{F^{p}H^{2p-1}(U,\partial U,\C )},\]
and hence

\begin{equation}
J^{p,n}(X) \simeq \frac {H^{2p-1} (U, \partial U,\C)}{F^p 
H^{2p-1}(U,\partial U,\C) + H^{2p-n-1} (X,\Z(p))} .\\
\end{equation}
Note that the formula in (6.4) generalizes the formula in (3.3) for 
the case
\( n=0 \) (\( \s _{X}^{0}=X \), \( \partial \s _{X}^{0}=0 \)). From 
(6.1),
we have a sequence analogous to (3.1), namely

\begin{equation}
\end{equation}
\[
0\to H^{2p-n-1}(X,\Z (p))\begin{array}[b]{c}
_{\tilde{\beta }}\\
\to \end{array}H^{2p-1}(U\m \partial U,\Z (p))\to \ker (\beta )\to 
H^{2p-n}(X,\Z (p)),\]
and modulo the image of \( \tilde{\beta } \), the higher Chow cycle 
\( \z \in CH^{p}_{\text {hom}}(X,n) \)
defines a class \( \Psi _{p,n}(\z )\in J^{p,n}(X) \), with the help 
of (6.4).

\section{\textbf{Comparing definitions}}

In this section it is proved that the constructions of \( AJ \) via 
explicit currents and via extension classes (in \( \S \S 5 \) and \( 
6 \), resp.) agree,
namely \begin{thm}
For \( X \) projective and \( [ \z ] \in CH^p_{hom}(X,n) \), \( \Psi
_{p,n}(\z )=\Phi _{p,n}(\z ). \)  If \( n\geq 1 \), then without
the "hom" assumption one has this equality modulo torsion (i.e., replacing
\( \Z \) by \( \Q \) in the target groups).  \end{thm}
\begin{proof}
We only need to show that the extension definition gives the same as
our explicit formula, since the equality between Bloch's map and the 
extension definition was already shown by Scholl \cite{Sch}.
We begin by picking apart
the last section's construction in some detail, and show in
particular that this construction agrees with the geometric
interpretation of the \(AJ\) map given in 5.8 above.

Recall \( \s ^{n}=(\P ^{1}\m \{1\}) \), \( \d \s 
^{n}=Alt_{n}(\{0,\infty \}\times (\P ^{1}\m \{1\})^{n-1}) \),
\( \s _{X}^{n}=X\times \s ^{n} \), etc. Let \( [\z ] \) be \( any \) 
class
\( \in CH^{p}(X,n) \). By normalization \cite{E-Z} we may choose \( 
\z \) to be a relative
cycle, so that all face intersections \( \z \cdot \partial \s 
^{n}_{X} \) are
zero as cycles; moreover, \( H^{2p}(\s ^{n}_{X}) \) injects into \( 
H^{2p}(\partial \s ^{n}_{X}) \)
(e.g. see diagram \( (6.2) \)). So the fundamental class of \( \z \) 
in \( H^{2p}_{|\z |}(\s _{X}^{n}) \)
goes to zero in \( H^{2p}(\s _{X}^{n}) \) and \( H^{2p}_{|\d \z |}(\d 
\s _{X}^{n}) \);
this accounts for \( \z \)'s determining a class in \( \ker (\beta ) 
\),
where
\[
\beta :\, H_{|\z |}^{2p}(\s ^{n}_{X})^{\circ }\to H^{2p}_{|\d \z 
|}(\d \s ^{n}_{X})^{\circ }.\]
None of this has anything to do with \( [\z ]\in CH_{hom}^{p}(X,n). \)

Now we describe the map \[
\ker (\beta )\to H^{2p-n}(X);\]
triviality of \( this \) \( map \) is what the ``hom'' indicates. 
Roughly speaking, this can be seen as also saying that the 
fundamental class
of \( \z \) in \( H^{2p}(\s ^{n}_{X},\d \s ^{n}_{X}) \) is trivial. 
This
corresponds to triviality of \( \z \) as a topological cycle in \( 
H_{2m+2n-2p}(X\times (\C ^{*},\{1\})^{n}), \)
which is best expressed by casting \( \z \) as a \( limit \) of 
topological
cycles \( \z _{\e } \) which \( are \) (modulo \( X\times \I ^{n}) \) 
boundaries
of topological \( (2m+2n-2p+1) \)-chains compactly supported on \( 
X\times (\C ^{*})^{n} \).

This map was defined via the following ``serpent'' of maps and lifts

\[
\begin{array}{ccccc}
& & H^{2p-1}(U) & \begin{array}{c}\longrightarrow\\ {<- - 
-}\end{array}
& H^{2p}_{|\z|}(\s^n_X)^{\circ} \\
&\\
& & \downarrow \\
&\\
H^{2p-1}(\d \s^n_X) & \begin{array}{c}\longrightarrow\\ {<- - 
-}\end{array} & H^{2p-1}(\d U) \\
&\\
\downarrow \\
&\\
H^{2p-n}(X) \\
\end{array}
\]
where the last (vertical) map takes \[
H^{2p-1}(\d \s ^{n}_{X})\cong H^{0}(\d \s ^{n})\otimes H^{2p-1}(X)\, 
\oplus \, H^{n-1}(\d \s ^{n})\otimes H^{2p-n}(X)\twoheadrightarrow 
H^{2p-n}(X),\]
which is to say \( \alpha \mapsto \pi _{*}(\alpha \v \di [\Omega 
^{n}]) \)
or (for topological cycles) \( \c \to \pi _{*}(\c \cdot \d 
T^{n}_{X}). \)\footnote{%
The point is that \( H_{|\d \s ^{n}|}^{n-1}((\P ^{1})^{n})\cong 
\left\{ H^{n-1}(\d \s ^{n})\right\} ^{\vee } \)
is generated by \( \d [\Omega ^{n}] \) (or Poincar\'e-dually by \( \d 
T^{n}) \);
therefore wedging with this and pushing down (integrating fiberwise) 
``removes''
the \( H^{n-1}(\d \s ^{n}) \) part from \( \alpha \).
} This latter description is more apparent from the equivalent 
homological serpent

\[
\begin{array}{ccccccc}
&&\begin{array}{c}{}_{\ker:}\\ {}^{\left[ H_{2m+2n-2p}(|\bar{\z}|) 
\to H_{2m+2n-2p}\left( (\P^1)^n_X , \I^n_X \right) 
\right]}\end{array}\\
&\\
&&\begin{array}{cc} \uparrow&|\\ |&(1)\\ |&\downarrow\end{array}\\
&\\
& & H_{2m+2m-2p+1} \left((\P^1)^n_X,|\bar{\z}| \cup \I^n_X \right)\\
&\\
& & \downarrow {2} \\
&\\
H_{2m+2n-2p-1} \left(\d (\P^1)^n_X , \d \bar{\I}^n_X \right) & 
\begin{array}{c}\to\\ ^{\leftarrow (3) -}\end{array} &
H_{2m+2n-2p-1} \left( \d (\P^1)^n_X , |\d \bar{\z}| \cup \d 
\bar{\I}^n_X \right) \\
&\\
\downarrow {4} \\
&\\
H_{2m+n-2p}(X). \end{array} \]

Tracing through, we have (since \( \z \) is in the upper right-hand 
kernel)
\[
\z =\d \Y \begin{array}[b]{c}
_{1}\\
\longrightarrow \end{array}\Y \begin{array}[b]{c}
_{2}\\
\longrightarrow \end{array}\begin{array}[t]{c}
\Y \cdot \d \s _{X}^{n}\\
^{\text {cycle\, mod\, }|\d \bar{\z }|\cup \d \bar{\I }_{X}^{n}}
\end{array}\begin{array}[b]{c}
_{3}\\
\longrightarrow \end{array}
\]
\[\begin{array}[t]{c}
\Y \cdot \d \s ^{n}_{X}\\
^{\text {cycle\, mod\, }\d \bar{\I }^{n}_{X}}
\end{array}\begin{array}[b]{c}
_{4}\\
\longrightarrow \end{array}\pi _{*}((\Y \cdot \d \s ^{n}_{X})\cdot \d 
T^{n}_{X})=\pi _{*}(\Y \cdot \d T^{n}_{X}).\]
Recalling that \( T_{\z }=\pi _{*}(\z \cdot T_{X}^{n}) \), (where 
$T^{n}_{X} = X\times T^{n}$), since \[
\d \pi _{*}(\Y \cdot T^{n}_{X})=\pi _{*}(\d \Y \cdot T^{n}_{X})\pm 
\pi _{*}(\Y \cdot \d T_{X}^{n})\]
we see that the image of \( \z \) under the above map is (mod 
coboundary
and up to sign) equivalent to \( T_{\z } \). Therefore we will 
write\footnote{%
or more formally, \( (2\pi \j)^{p}T:\, \ker (\beta )\to H^{2p-n}(X,\Z 
(p)) \).
} \[
T:\, \ker (\beta )\to H^{2p-n}(X).\]

If \( [\z ]\in CH^{p}_{hom}(X,n) \) then \( T_{\z }\sim 0 \).
Thus we get an element in the right-hand term of the twisted sequence 
{[}from (6.5){]}
\[
0\to H^{2p-n-1}(X,\Z (p))\begin{array}[b]{c}
A\\
\longrightarrow \end{array}H^{2p-1}(U,\d U)_{\Z(p)}\begin{array}[b]{c}
B\\
\longrightarrow \end{array}\ker (T)_{\Z (p)} \to 0;\]
we show how to lift it to the center term. The map \( B \) is the 
composite
\[
H^{2p-1}(U,\d U)\to H^{2p-1}(U)\to H^{2p}_{|\z |}(\s ^{n}_{X})\]
or (homologically)
\[
\begin{array}{cc} &H_{2m+2n-2p}(|\bar{\z }|)\\
&\\
&\begin{array}{cc} \uparrow&|\\ \d&(1)\\
|&\downarrow\end{array}\\
&\\
H_{2m+2n-2p+1}(X\times (\C ^{*})^{n},|\z |\cup \I ^{n}_{X}) 
\begin{array}[t]{c}
\longrightarrow \\
^{\leftarrow(2)-}
\end{array}&H_{2m+2n-2p+1}(X\times (\P ^{1})^{n},|\bar{\z }|\cup 
\bar{\I }_{X}^{n})\end{array}
\]
where the lifts indicated (possible because of the exact sequence) 
take \( \z \)
to a bounding chain \( \Gamma \) (where \( \d \Gamma =\z \) mod \( \I 
^{n}_{X} \))
and then to a limit of bounding chains \( \Gamma _{\e } \) compactly 
supported
on \( X\times (\C ^{*})^{n} \) (with \( \d \Gamma _{\e }=\z _{\e } \) 
mod \( \I ^{n}_{X} \)) as described in \( \S 5 \).

We will write \( \zeta _{\z } \) for the \( image \) of the integral 
lift
\[ \lim _{\e \to 0}\Gamma _{\e } \ \text{in} \ H^{2p-1}(U,\d U)_{\C 
}/F^{p}H^{2p-1}(U,\d U). \]
Its preimage \( \tilde{\beta }^{-1}(\zeta _{\z }) \) under
\[
\tilde{\beta }:\, \frac{H^{2p-n-1}(X)_{\C 
}}{F^{p}H^{2p-n-1}(X)}\begin{array}[b]{c}
_{\cong }\\
\longrightarrow \end{array}\frac{H^{2p-1}(U,\d U)}{F^{p}H^{2p-1}(U,\d 
U)},\]
taken modulo \( H^{2p-n-1}(X,\Z (p)) \), gives \( \Psi _{p,n}(\z ). 
\) Now
\( \tilde{\beta } \) dualizes to
\[
\tilde{\beta }^{\vee }:\, F^{n+m-p+1}H^{2n+2m-2p+1}\left( X\times (\C 
^{*})^{n},\, |W|\cup \I ^{n}_{X}\right)_{\C }\begin{array}[b]{c}
_{(\cong )}\\
\longrightarrow \end{array} \hskip1.5in
\]
\[
\hskip2in F^{m-p+1}H^{n+2m-2p+1}(X)_{\C }.\]
If we think of \( \tilde{\beta }^{_{-1}}(\zeta _{\z })\in \left\{ 
F^{m-p+1}H^{2m-2p+n+1}(X)_{\C }\right\} ^{\vee } \)
as a functional on forms, then it (and hence \( \Psi _{p,n}(\z ) \) ) is
computed on \( \omega \) by
\[
(*)\, \, \, \, \, \, \, \, \, \, \, \, \, \, \, \, [\tilde{\beta 
}^{-1}(\zeta _{\z })]\omega \, =\, \zeta _{\z }[(\tilde{\beta }^{\vee 
})^{-1}\omega ].\]
(This essentially comes from Carlson's theory \cite{C}.) 

It remains to trace through \( \tilde{\beta }^{\vee } \). As \( 
\tilde{\beta } \)
is given by the composition

\[
\begin{array}{ccccc}
H^{2p-n-1})(X) & \hookrightarrow & H^{2p-2}(\d \s^n_X) & \to & 
H^{2p-1}(U,\d U) \\
&\\
{\buildrel{{}_{\cap}\hskip.045in }\over\downarrow} & \nearrow {\cong} 
& &\searrow & \uparrow \\
&\\
H^{2p-n-1}(X)\oplus H^{2p-2}(X) & & & & H^{2p-2}(\d U) \\
\end{array} \]
\( \tilde{\beta }^{\vee } \) must be
\[
H^{2m-2p+n+1}(X)\twoheadleftarrow H_{|\d \s ^{n}_{X}|}(X\times (\P 
^{1})^{n})\begin{array}[b]{c}
_{\di [\cdot ]}\\
\longleftarrow \end{array}H^{2m+2n-2p+1}(X\times (\C ^{*})^{n},\, |\z 
|\cup \I ^{n}_{X})\]
where the second map sends \( \pi ^{*}\omega \v \di [\Omega ^{n}]\, 
\mapsto \, \omega \),
and the group on the r.h.s. we take to be represented by forms 
pulling back
to \( 0 \) along \( |\z | \). If \( \omega \in \Gamma 
(F^{m-p+1}\Omega ^{2m-2p+n+1}_{X^{\infty }}) \)
is (\( \di \)-)closed then \( \pi ^{*}\omega \v \Omega ^{n} \) is of 
type
\( F^{m+n-p+1} \) (whereas \( \dim _{\C }\z =m+n-p \)), and so gives 
a lift
of \( \omega \) to the r.h.s.; therefore we write \( (\tilde{\beta 
}^{\vee })^{-1}\omega =\pi ^{*}\omega \v \Omega ^{n} \).
So \( (*) \) is just \[
\lim _{\e \to 0}\int _{\Gamma _{\e }}\pi _{1}^{*}\omega \v \pi 
^{*}_{2}\Omega ^{n},\]
identifying \( \Psi _{p,n}(\omega ) \) with the \( AJ \) for relative 
cycles
\( \in Z^{p}(\s ^{n}_{X},\d \s ^{n}_{X}) \) as described in \( \S 5 
\). We
already know this equates with \( \Phi _{p,n}(\z ) \) (or \( AJ(\z ) 
\)) by Proposition~\ref{currentequate}. This completes the proof.
\end{proof}
\begin{rem*}
A few words need to be said regarding the ``\( \lim _{\e \to 0} \)''. 
The
``relative quasiprojective variety'' \[
(**)\, \, \, \, \, \, \, \, (X\times (\C ^{*})^{n},\, |\z |\cup \I 
^{n}_{X})\]
is Poincar\'e-dual to \( (U,\d U) \); here these two play roles 
analogous to
those played, respectively, by \( (X,\, |\z |) \) and \( (X\m |\z |) 
\) in
\( \S 3 \). We need to be able to pair forms on \( (**) \) {[}which 
pull back
to \( 0 \) along \( |\z |\cup \I ^{n}_{X} \) but have poles along \( 
\d \s ^{n}_{X} \){]}
with topological cycles there; this is the reason for using limits of 
chains
to compute homology (the integrals only make sense as a limit). 
Furthermore
this needs to be done in such a way that coboundaries \& topological 
cycles
{[}resp. cocycles \& topological boundaries{]} pair to zero; 
referring to \( \S 5 \),
the fact that \( \lim _{\e \to 0}\int _{\W _{\e }}\pi _{\s 
}^{*}\left\{ \begin{array}{c}
\Omega ^{n}\\
R^{n}
\end{array}\right\} \v \pi ^{*}_{X}\alpha \, =\, 0 \) for any \( 
C^{\infty } \)-form \( \alpha \) ensures that this condition
is met (as far as needed for the arguments here).
\end{rem*}

\end{document}